\documentclass[11pt,final]{siamltex}%
\usepackage{amsmath}
\usepackage{amsfonts}
\usepackage{amssymb}
\usepackage{graphicx}%
\providecommand{\U}[1]{\protect\rule{.1in}{.1in}}
\providecommand{\U}[1]{\protect\rule{.1in}{.1in}}
\begin{document}

\title{ A parallel method for solving Laplace equations with Dirichlet data using
local boundary integral equations and random walks \thanks{Submitted to SIAM
J. Scientific Computing in revision on October 4, 2012. } }
\author{Chanhao Yan\thanks{State Key Lab. of ASIC \& System Fudan Univ. Shanghai,
China (\texttt{yanch@fudan.edu.cn})}
\and Wei Cai \thanks{Dept. of Mathematics \& Statistics University of North
Carolina at Charlotte, USA}. (\texttt{wcai@uncc.edu})
\and Xuan Zeng\thanks{State Key Lab. of ASIC \& System Fudan Univ. Shanghai, China
(\texttt{xzeng@fudan.edu.cn})} }
\maketitle

\begin{abstract}
In this paper, we will present a new approach for solving Laplace equations in
general 3-D domains. The approach is based on a local computation method for
the DtN mapping of the Laplace equation by combining a deterministic (local)
boundary integral equation method and the probabilistic Feynman-Kac formula of
PDE solutions. This hybridization produces a parallel algorithm where the bulk
of the computation has no need for data communications. Given the Dirichlet
data of the solution on a domain boundary, a local boundary integral equation
(BIE) will be established over the boundary of a local region formed by a
hemisphere superimposed on the domain boundary. By using a homogeneous
Dirichlet Green's function for the whole sphere, the resulting BIE will
involve only Dirichlet data (solution value) over the hemisphere surface, but
over the patch of the domain boundary intersected by the hemisphere, both
Dirichlet and Neumann data will be used. Then, firstly, the solution value on
the hemisphere surface is computed by the Feynman-Kac formula, which will be
implemented by a Monte Carlo walk on spheres (WOS) algorithm. Secondly, a
boundary collocation method is applied to solve the integral equation on the
aforementioned local patch of the domain boundary to yield the required
Neumann data there. As a result, a local method of finding the DtN mapping is
obtained, which can be used to find all the Neumann data on the whole domain
boundary in a parallel manner. Finally, the potential solution in the whole
space can be computed by an integral representation using both the Dirichlet
and Neumann data over the domain boundary.

\end{abstract}



\section{Introduction}

Fast and parallel scalable solvers for 3-D Poisson and modified Helmholtz
partial differential equations (PDEs) constitute the major computational cost
for many large-scale scientific computational problems, such as
Poisson/Helmholtz solvers in projection type methods of incompressible flows
\cite{chorin68}\cite{temam84}, electrostatic potential problems in molecular
biology, and enforcing divergence-free constraints of magnetic fields in
plasma MHD simulations, etc. For the electrostatic capacitance problems for
conductors, boundary element methods (BEMs) or finite element methods (FEMs)
are often used to compute the charge density in the engineering community, for
example, the indirect BEM FastCap \cite{white91}\cite{white94}, the direct BEM
QMM-BEM \cite{yu04}, hierarchical extractors HiCap and PhiCap \cite{shi02}%
\cite{shi05}, and the parallel adaptive FEM ParAFEMCap \cite{zhu12}, etc. BEMs
\cite{brebbia78} need to discretize entire conductor surfaces, sometimes even
the dielectric interfaces, into small panels, and construct a linear system by
the method of moments or collocation methods. These deterministic methods are
highly accurate and versatile, but are global, i.e., even if the charge
density at only one point is required, a full linear system has to be
constructed and solved. In general, the resulting linear algebraic systems are
solved by iterative methods such as the multi-grid methods \cite{brandt82} or
the domain decomposition methods \cite{widlund05}, either as a solver or as a
pre-conditioner. Meanwhile, for integral equation discretization, the fast
multipole method (FMM) \cite{Greengard} can be used in conjunction with a
Krylov subspace iterative solver. All these solvers are $O(N)$ in principle
and iterative in nature, and require expensive surface or volume meshes. The
parallel scalability of these solvers on a large number of processors poses
many challenges and is the subject of intensive research.

In contrast, random methods can give local solutions of PDEs \cite{given02}%
\cite{shebalfeld}\cite{hwang06}, and they have been applied to obtain
solutions at specific sites for many real world problems such as modern VLSI
chips with millions of circuit elements in the area of chip design industry.
\ For instance, the QuickCap, as the chip industry's gold standards produced
by the leading EDA companies Synopsys, is a random method. The key advantage
of the random methods is their localization. For example, QuickCap
\cite{coz92}\cite{coz98} can calculate the potential or charge density at only
one point locally without finding the solution elsewhere. Usually, random
methods are based on the Feynman-Kac probabilistic formula and the potential
(or charge density) is expressed as a weighted average of the boundary values
\cite{hwang06}. The Feynman-Kac formula allows a local solution of the PDE,
and fast sampling techniques of the diffusion paths with the walk on sphere
(WOS) methods are available for simple PDEs such as Laplace or modified
Helmholtz equations. However, it is impractical to use the probabilistic
formula to find the solution of these PDEs in the whole space as too much
sampling will be needed. \

\bigskip For current multi-core petaflops per second computing platforms, the
scalability of the algorithms becomes the major concern for the development of
new algorithms. Much research has been done in order to achieve such a
scalability and parallelism in the above deterministic algorithms for
simulation capability for realistic engineering and scientific problems. To
meet this challenge, in this paper, we will propose a hybrid method for
computing the Neumann data (derivative) of the solution from its Dirichlet
data by combining the probabilistic Feynman-Kac formula and a deterministic
local integral equation over a domain boundary $\partial\Omega$. The hybrid
method will allow us to get the Neumann data efficiently over a local patch of
the domain boundary, which will result in a simple intrinsic parallel method
for solving the complete potential problems in general 3-D domains through an
integral representation of the available Dirichlet and Neumann data.

The rest of the paper will be organized in the following sections. Firstly, we
will present some background material on the Dirichlet to Neumann (DtN)
mapping and also the Feynman-Kac probabilistic solution of elliptic PDEs.
Secondly, we will review a related last--passage random walk method proposed
in \cite{given02} which calculates the Neumann data (charge distribution)
\textit{at one single point} over a flat surface where the Dirichlet data is a
constant. Even though this is a very limited case for the DtN problem, it
demonstrates some key issues and difficulties in how to use the Feynman-Kac
formula and the WOS in finding the Neumann data. Thirdly, we will present our
hybrid method, which allows the calculation of the Neumann data for a general
Dirichlet data on the flat surface. Then, in Section 4 the hybrid method is
extended to calculate the Neumann data over \textit{a patch} of the boundary
for arbitrary Dirichlet data and curved boundaries. In Section 5 numerical
tests will be presented to show the accuracy and potential of the proposed
method. Finally, conclusions and discussions for open research issues and
parallel aspect of the proposed method will be given in Section 6.

\section{Background on DtN mapping and solutions of potential equations}

The DtN mapping between the Dirichlet data (solution value) and the Neumann
data (the normal derivative of the solution) of a Poisson equation is relevant
in both engineering applications and mathematical study of elliptic PDEs. In
the electrostatic potential problems, the surface charge distribution
$\sigma_{s}$ on the surface $\partial\Omega$ of a conductor $\Omega$, as
required in the capacitance calculation of conductive interconnects in VLSI
chips, is exactly the normal derivative of the electrostatic potential $u$ as
implied from Gauss's law for the electric field $\mathbf{E=-\nabla}u$, i.e.,%

\begin{equation}
\sigma_{s}=\mathbf{E\cdot n|}_{\partial\Omega}\mathbf{=-}\left.
\frac{\partial u}{\partial\mathbf{n}}\right\vert _{\partial\Omega}.
\label{bw00}%
\end{equation}


On the other hand, the DtN mapping also plays an important role in the study
of the Poisson equations.\ As the inhomogeneous right--hand--side of a Poisson
equation is usually known, we could use a simple subtraction technique to
reduce the Poisson equation to a Laplace equation, but with a modified
boundary data. Therefore, in the rest of this paper we will present our method
for the Laplace equation in a domain $\Omega$ where a general Dirichlet data
is given on the boundary $\partial\Omega$. If we are able to compute the
Neumann data for the given Dirichlet data, namely the following DtN mapping:%

\begin{equation}
\text{DtN: \ \ \ \ }u|_{\partial\Omega}\rightarrow\left.  \frac{\partial
u}{\partial\mathbf{n}}\right\vert _{\partial\Omega}, \label{bw01}%
\end{equation}
then, the solution $u(\mathbf{x})$ at any point $\mathbf{x}$ in the whole
space can be found simply by the following integral representation:%

\begin{equation}
u(\mathbf{x})=%
{\displaystyle\int\nolimits_{\partial\Omega}}
G(\mathbf{x},\mathbf{y})\frac{\partial u(\mathbf{y})}{\partial\mathbf{n}%
_{\mathbf{y}}}ds_{y}-%
{\displaystyle\int\nolimits_{\partial\Omega}}
\frac{\partial G(\mathbf{x,y})}{\partial\mathbf{n}_{\mathbf{y}}}%
u(\mathbf{y})ds_{y},\text{ \ \ }\mathbf{x}\in\mathbb{R}^{3}\backslash
\partial\Omega, \label{bw03}%
\end{equation}
where $G(\mathbf{x},\mathbf{y})$ is the fundamental solution to the Laplace
operator, namely,%

\begin{equation}
G(\mathbf{x},\mathbf{y})=\frac{1}{4\pi}\frac{1}{|\mathbf{x-y}|}. \label{bw05}%
\end{equation}

A similar NtD (Neumann to Dirichlet) mapping from the Neumann data to the
Dirichlet data can also be defined if the Neumann data yields a unique
solution to the PDE. In either case, with both Dirichlet and Neumann data at
hand, the solution of a Laplace equation can be obtained by the representation
formula in ($\ref{bw03}$).

Therefore, by finding the DtN or NtD mapping of the relevant elliptic PDE
solutions in an efficient manner, we could produce fast numerical methods for
many electrical engineering and fluid mechanics applications.

The Feynman-Kac formula \cite{Freidlin85}\cite{friedman06} relates the Ito
diffusion paths to the solution $u(\mathbf{x})$ of the following general
elliptic problem%
\begin{align}
L(u)  &  \equiv%
{\displaystyle\sum\limits_{i=1}^{3}}
b_{i}(\mathbf{x})\frac{\partial u}{\partial x_{i}}+%
{\displaystyle\sum\limits_{i,j=1}^{3}}
a_{ij}(\mathbf{x})\frac{\partial^{2}u}{\partial x_{i}\partial x_{j}%
}=f(\mathbf{x}),\text{ \ \ \ }\mathbf{x}\in\Omega,\nonumber\\
u|_{\partial\Omega}  &  =\phi(\mathbf{x}),\text{ \ }\mathbf{x}\in
\partial\Omega, \label{bw07}%
\end{align}
where $L$ is an uniformly elliptic differential operator, i.e.%

\begin{equation}%
{\displaystyle\sum\limits_{i,j=1}^{3}}
a_{ij}(\mathbf{x})\xi_{i}\xi_{j}\geq\mu|\boldsymbol{\xi }|^{2},\text{ \ if
\ \ \ }\mathbf{x\in}\Omega,\boldsymbol{\xi }\in R^{3}\text{ \ \ }(\mu>0),
\end{equation}
and $a_{ij}(\mathbf{x})$ and $b_{i}(\mathbf{x})$ are uniformly Lipschitz
continuous in $\overline{\Omega}.$ Also the domain $\Omega$ is assumed to have
a $C^{2}$ boundary and the boundary data $\phi\in$ $C^{0}$ $(\partial\Omega)$.

If $X_{t}(\omega)$ is an Ito diffusion defined by the following stochastic
differential equation%

\begin{equation}
dX_{t}=b(X_{t})dt+\alpha(X_{t})dB_{t}, \label{bw09}%
\end{equation}
where $B_{t}(\omega)$ is the Brownian motion, $[a_{ij}]=\frac{1}{2}%
\alpha(x)\alpha^{\mathrm{T}}(x)$, then, the following Feynman-Kac formula
gives a probabilistic solution for (\ref{bw07}) as%

\begin{equation}
u(\mathbf{x})=E^{x}\left(  \phi(X_{\tau_{\Omega}})\right)  +E^{x}\left[
{\displaystyle\int\nolimits_{0}^{\tau_{\Omega}}}
f(X_{t})dt\right]  ,\label{bw11}%
\end{equation}
where the expectation is taken over all sample paths $X_{t=0}(\omega
)=\mathbf{x}$ and $\tau_{\Omega}$ is the first hit time (or exit time) of the
domain $\Omega$. For the purpose of this paper, we will only consider
(\ref{bw11}) for Laplace equations ($f=0$).

For the Laplace equation, the Ito diffusion is just the Brownian motion and
the solution can be simply rewritten in terms of a harmonic measure, which
measures the probability of the Brownian paths hitting a given area on the
boundary surface,%

\begin{equation}
u(\mathbf{x})=E^{x}(\phi(X_{\tau_{\Omega}}))=%
{\displaystyle\int\nolimits_{\partial\Omega}}
\phi(\mathbf{y})d\mu_{\Omega}^{x}, \label{bw13}%
\end{equation}
where%

\begin{equation}
\mu_{\Omega}^{x}(F)=P^{x}\{\omega|X_{\tau_{\Omega}}(\omega)\in F,X_{0}%
(\omega)=\mathbf{x}\},\ F\subset\partial\Omega,\text{ \ }\mathbf{x}\in\Omega.
\label{bw15}%
\end{equation}
\qquad\qquad

The harmonic measure can be shown to be related to Green's function
$g(\mathbf{x},\mathbf{y})$ of the Laplace equation in the domain $\Omega$ with
a homogeneous boundary condition, i.e.,%
\begin{align}
-\Delta g(\mathbf{x},\mathbf{y})  &  =\delta(\mathbf{x}-\mathbf{y}),\text{
\ }\mathbf{x}\in\Omega,\nonumber\\
g(\mathbf{x},\mathbf{y})|_{\mathbf{x}\in\partial\Omega}  &  =0. \label{bw17}%
\end{align}

By comparing ($\ref{bw13}$) with the following integral representation of the
solution of the Laplace equation with Green's function $g(\mathbf{x}%
,\mathbf{y})$,%

\begin{equation}
u(\mathbf{x})=-%
{\displaystyle\int\nolimits_{\partial\Omega}}
\phi(\mathbf{y})\frac{\partial g(\mathbf{x},\mathbf{y})}{\partial
\mathbf{n}_{\mathbf{y}}}ds_{y},\label{bw19}%
\end{equation}
we conclude that the hitting probability, now denoted as $p(\mathbf{x,y}%
)ds_{y}=\mu_{\Omega}^{x}([\mathbf{y},\mathbf{y}+ds_{y}])$, has the following
connection to Green's function of the domain $\Omega$ \cite{Chung95},%

\begin{equation}
p(\mathbf{x,y})=-\frac{\partial g(\mathbf{x},\mathbf{y})}{\partial
\mathbf{n}_{\mathbf{y}}}. \label{bw21}%
\end{equation}

\begin{figure}[ptb]
\begin{center}
\includegraphics[width=2.8in ]{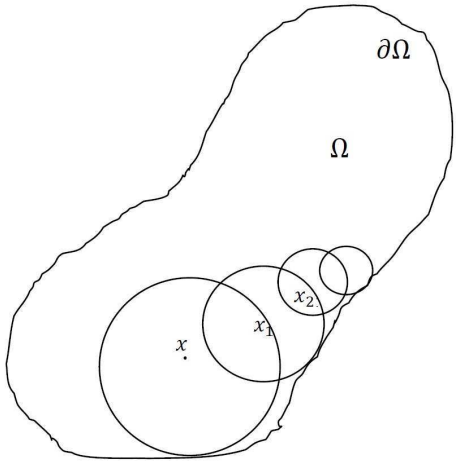}
\end{center}
\caption{WOS sampling of Brownian paths}%
\label{fig_wos}%
\end{figure}

Therefore, if the domain is a ball centered at $\mathbf{x}$ where a path
starts, we have a uniform probability for the path to hit the surface of the
ball. This fact will be a key factor in the design of random walk on spheres
(WOS), which allows us to describe the Brownian motion and its exit location
without explicitly finding its trajectory. Instead, a sequence of walks or
jumps over spheres will allow the Brownian path to hit the boundary
$\partial\Omega$ (for practical purpose, within an absorption $\varepsilon
$-shell as proposed in \cite{Mascag03}). Specifically, as indicated by
(\ref{bw21}), the probability of a Brownian path hitting on a spherical
surface is given exactly by the normal derivative of Green's function of the
sphere (with a homogeneous boundary condition). Therefore, if we draw a ball
centered at the starting point $\mathbf{x}$ of a Brownian path, it will hit
the spherical surface with a uniform probability as long as the ball does not
intersect with the domain boundary $\partial\Omega$. So, we can make a jump
for the Brownian particle to $\mathbf{x}_{1}$, sampled with a uniform
distribution on the spherical surface. Next, a second ball now centered at
$\mathbf{x}_{1}$ will be drawn, not intersecting with the domain boundary
$\partial\Omega$, and the Brownian particle can make a second jump to
$\mathbf{x}_{2}$ on the surface of the second ball. This procedure (as
illustrated in Fig. \ref{fig_wos}, termed as WOS) \cite{Muller56}%
\cite{given02}\cite{mascagni04} is repeated until the Brownian particle hits
the \text{boundary of }$\Omega$ (within the $\varepsilon$-shell of absorption)
where it is denoted as $\mathbf{x}_{\tau_{\Omega}}$ and the value of the
boundary data $\phi(\mathbf{x}_{\tau_{\Omega}})$ will be recorded and
eventually all such data will be used to compute the expectation in
(\ref{bw13}). In real applications, due to the relation between Green's
function $g(\mathbf{x,y})$ of a domain and the hitting probability, Green's
Function First Passage (GFFP) methods for shapes other than spheres such as
rectangles in softwares including QuickCap \cite{coz92}\cite{coz98} have been
used to find capacitances of conductors in interconnect layouts, which are
generally of rectangular shapes.

Moreover, in applying the Feynman-Kac formula (\ref{bw13}) to find the
potential in the exterior domain of infinite extent (with a vanishing
condition for the potential at the infinity), as some paths will go to
infinity, a truncation procedure by a large sphere is used in our simulation
of WOS where trajectories outside the big sphere will be ignored and
considered as ending at infinity where the potential value vanishes.
Theoretical estimate on the size of the truncation sphere can be found in
\cite{shebalfeld95}.

\section{Finding the Neumann data at one point on a flat boundary}

\subsection{Review of the last-passage algorithm for charge density}

In this subsection, we will review the last-passage Monte Carlo algorithm
proposed in \cite{given02} for charge density, namely the Neumann data, at one
point on a flat \textit{conducting} surface.

\begin{figure}[ptb]
\begin{center}
\includegraphics[width=4in ]{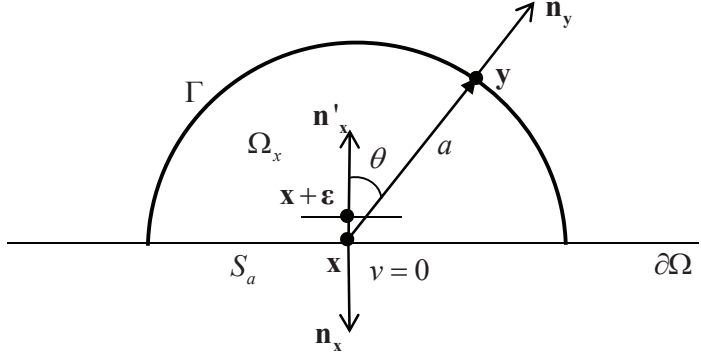}
\end{center}
\caption{Last--passage for finding the Neumann data at one point}%
\label{fig_lastpass}%
\end{figure}

For a flat portion of the boundary $\partial\Omega$ of a domain $\Omega
=\{z<0\}$ in the 3-D space held at a constant potential, we like to compute
the charge density at a point $\mathbf{x\in}$ $\partial\Omega$, namely, the
normal derivation of the exterior potential outside $\Omega$. In the
last-passage method, a hemisphere is constructed with a radius $a$ centered at
$\mathbf{x}$ as shown in Fig. \ref{fig_lastpass}. The upper hemispherical
surface outside $\Omega$ is denoted as $\Gamma$ and the 2-D disk of radius $a$
centered at $\mathbf{x}$ from the intersection of the hemisphere and the
conducting boundary $\partial\Omega$ is denoted as%
\begin{equation}
S_{a}\equiv S_{a}(\mathbf{x}). \label{bw20}%
\end{equation}

In Section 2, the equivalence between the electrostatic potential $\ u$ (which
is assumed at value $1$ on the conductor surface $\partial\Omega$) and
diffusion problems is given. In the last-passage method \cite{given02}, the
quantity $v(\mathbf{x})\equiv1-u(\mathbf{x})$ is considered, which will
satisfy%
\begin{equation}
v(\mathbf{x})=0,\text{ \ }\mathbf{x}\in S_{a}\label{bw22}%
\end{equation}
and $v=1$ at infinity (or on an infinitely large sphere). By viewing
$v(\mathbf{x}+\varepsilon)$ as the probability of a Brownian particle near the
conducting surface $\partial\Omega$ starting at $\mathbf{x}+\varepsilon$
diffusing to infinity without ever coming back to the conducting surface, it
was shown in  \cite{given02} that the following probabilistic expression for
$v(\mathbf{x}+\varepsilon)$ holds:%

\begin{equation}
v(\mathbf{x}+\varepsilon)\equiv1-u(\mathbf{x}+\varepsilon)=-%
{\displaystyle\int\nolimits_{\Gamma}}
\widehat{g}(\mathbf{x}+\varepsilon,\mathbf{y})p_{y\infty}ds_{y}, \label{bw23}%
\end{equation}
where $p_{y\infty}$ is the probability of a Brownian particle starting at
$\mathbf{y}$ and diffusing to infinity without ever coming back to the
conducting surface $\partial\Omega$; thus, $p_{y\infty}=0$ if $\mathbf{y}\in
S_{a}$. In (\ref{bw23}), the integral over $\Gamma$ expresses the Markov
property of the diffusing particles from $\mathbf{x}+\varepsilon$ to infinity
with an intermediate stop on $\Gamma$. Specifically, $\widehat{g}%
(\mathbf{x}+\varepsilon,\mathbf{y})$ gives the probability of a Brownian
particle starting from $\mathbf{x}+\varepsilon$ and hitting the boundary of
$\Gamma$, which is given by (\ref{bw21}) via a homogeneous Green's function
for the hemisphere over $S_{a}$, namely,$\ $%
\begin{equation}
\widehat{g}(\mathbf{x}+\varepsilon,y)=\frac{\partial g}{\partial\mathbf{n}%
_{y}}(\mathbf{x}+\varepsilon,\mathbf{y}), \label{bw25}%
\end{equation}
and $g(\mathbf{x}+\varepsilon,y)$ is defined in (\ref{bw17}) for the
hemisphere, whose analytical form can be obtained by an image method with
respect to the spherical surface first, then to the plane $z=0$, resulting in
the use of three images. Specifically, we have%
\begin{equation}
g(\mathbf{x},\mathbf{x}_{\mathrm{s}})=\frac{1}{4\pi}\frac{1}{|\mathbf{x}%
-\mathbf{x}_{s}|}+\frac{1}{4\pi}\frac{q_{\mathrm{k}}}{|\mathbf{x}%
-\mathbf{x}_{\mathrm{k}}|}+\frac{1}{4\pi}\frac{q_{\overline{\mathrm{s}}}%
}{|\mathbf{x}-\mathbf{x}_{\overline{\mathrm{s}}}|}+\frac{1}{4\pi}%
\frac{q_{\overline{\mathrm{k}}}}{|\mathbf{x}-\mathbf{x}_{\overline{\mathrm{k}%
}}|}, \label{bw26}%
\end{equation}
where in the spherical coordinates the source location is $\mathbf{x}%
_{\mathrm{s}}=(\rho_{\mathrm{s}},\theta_{\mathrm{s}},\phi_{\mathrm{s}}),$ the
Kelvin image location with respect to the sphere is $\mathbf{x}_{\mathrm{k}%
}=(a^{2}/\rho_{\mathrm{s}},\theta_{\mathrm{s}},\phi_{\mathrm{s}}),$ and their
mirror image locations with respect to the plane $z=0$ are $\mathbf{x}%
_{\overline{\mathrm{s}}}=(\rho_{\mathrm{s}},\pi-\theta_{\mathrm{s}}%
,\phi_{\mathrm{s}}),\mathbf{x}_{\overline{\mathrm{k}}}=(a^{2}/\rho
_{\mathrm{s}},\pi-\theta_{\mathrm{s}},\phi_{\mathrm{s}}),$ respectively.
Meanwhile, the corresponding charges are $q_{\mathrm{k}}=-a/\rho_{\mathrm{s}%
},$ $q_{\overline{\mathrm{s}}}=-1,$ and $q_{\overline{\mathrm{k}}}%
=a/\rho_{\mathrm{s}},$ respectively.

Now, to get the charge distribution $\sigma_{s}$ (normal derivative), we use
the relation in (\ref{bw00}), and we have%

\begin{equation}
\sigma_{s}=-\underset{\varepsilon\rightarrow0}{\lim}\mathbf{n}_{x+\varepsilon
}\cdot\mathbf{E(\mathbf{x}+\varepsilon)=}\ \underset{\varepsilon\rightarrow
0}{\lim}\frac{\partial u(\mathbf{x}+\varepsilon)}{\partial\mathbf{n}_{x}%
}=\frac{\partial u(\mathbf{x})}{\partial\mathbf{n}_{x}}, \label{bw27}%
\end{equation}
and%

\begin{equation}
\frac{\partial u(\mathbf{x})}{\partial\mathbf{n}_{x}}=%
{\displaystyle\int\nolimits_{\Gamma}}
h(\mathbf{x}+\varepsilon,\mathbf{y})p_{y\infty}ds_{y}\equiv\Sigma
_{\mathrm{LP}}, \label{bw29}%
\end{equation}
where the shorthand $\Sigma_{\mathrm{LP}}$ is introduced for the integral over
$\Gamma$ for latter use, and%

\begin{equation}
h(\mathbf{x,y})=\frac{\partial^{2}}{\partial\mathbf{n}_{x}\partial
\mathbf{n}_{y}}g(\mathbf{x}+\varepsilon,\mathbf{y}).\label{bw31}%
\end{equation}
The weight function $h(\mathbf{x,y})$ can be analytically computed for the hemisphere%

\begin{equation}
h(\mathbf{x,y})=\frac{3}{2\pi}\frac{\cos\theta}{a^{3}}, \label{bw33}%
\end{equation}
where $\theta$ is the angle between the two normal vectors $\mathbf{n}%
_{x}^{\prime}$ and $\mathbf{n}_{y}$ on the boundary $\Gamma$ as shown in Fig.
\ref{fig_lastpass}.

Next, we only need to compute $p_{y\infty}$ which is the probability of a
Brownian particle starting from $\mathbf{y}\in\Gamma$ and diffusing to
infinity without ever returning to the conductor surface $\partial\Omega.$ Due
to the homogeneity of the Brownian motion in the external domain $\Omega
^{c}=\{z>0\}$, the WOS in Section 2 can be used to calculate this probability.
The integral in (\ref{bw29}) over $\Gamma$ could be approximated by a Gauss
quadrature as both $h(\mathbf{x,y})$ and $p_{y\infty}$ can be considered as
smooth functions of $\mathbf{y}\in\Gamma.$ Nonetheless, in \cite{given02}, the
integral $\Sigma_{\mathrm{LP}}$ is computed by first distributing $N$
particles at locations over $\Gamma$ based on a distribution density derived
from (\ref{bw33}), and then starting a Brownian diffusion path from each of
those locations. The number of paths along which the particles will diffuse to
infinity (in practice, when it hits a very large ball) is recorded as
$N_{\inf}$. Then, we have the following estimate%

\begin{equation}
\Sigma_{\mathrm{LP}}\simeq\frac{3}{2a}\frac{N_{\inf}}{N}. \label{bw35}%
\end{equation}

The key equation in the last--passage algorithm is (\ref{bw29}), which is
based on (\ref{bw23}) provided that the potential solution $v(\mathbf{x}%
)=0,\mathbf{x}\in S_{a}$ on the conductor surface as indicated in
(\ref{bw22}). Therefore, for general non-constant Dirichlet boundary data, the
last--passage method will not be applicable. In fact, the charge density at
$\mathbf{x}$ will be influenced by the potential value on all domain boundaries.

\subsection{BIE-WOS Method: Combining a BIE and the Monte Carlo WOS method}

For the last-passage method discussed above, the algorithm (\ref{bw29}) is
obtained by the isomorphism between the electrostatic potential and diffusion
problems. The limitation of the last--passage method is that it is only
applicable to the situation of constant Dirichlet data and a flat boundary. In
this section, we will adopt a different approach based on a boundary integral
equation (BIE) representation of the charge density (the Neumann data) on the
surface at a given point using potential over a small hemisphere; the latter
can then be computed by the random WOS method. As a result, this new approach,
a hybrid method of deterministic and random approaches, will be able to handle
general variable Dirichlet boundary data, and later in Section 4 is also
extended to curved boundaries.

Let us denote by $\Omega_{x}$ the domain formed by the hemisphere of radius
$a$ centered at $\mathbf{x}$ over the flat boundary $S_{a}$ as in Fig.
\ref{fig_lastpass}. By applying the integral representation (\ref{bw03}) of
the Laplace equation with the afore mentioned Green's
function\ $g(\mathbf{x,y})$ in (\ref{bw26}) for the domain $\Omega_{x}$ with a
homogeneous Dirichlet boundary condition, due to the zero boundary value of
Green's function $g(\mathbf{x,y})$, we have%
\begin{equation}
u(\mathbf{x}^{\prime})=-%
{\displaystyle\int\nolimits_{\Gamma\cup S_{a}}}
\frac{\partial g(\mathbf{x}^{\prime}\mathbf{,y})}{\partial\mathbf{n}%
_{\mathbf{y}}}u(\mathbf{y})ds,\text{ \ \ }\mathbf{x}^{\prime}\in\Omega_{x},
\label{bw37}%
\end{equation}
where $\Gamma$ again is the surface of the upper hemisphere and $S_{a}$ is the
disk of radius $a$ centered at $\mathbf{x}$. In order to obtain the normal
derivative of $u$ at $\mathbf{x}$, we simply take the derivative with respect
to $\mathbf{x}^{\prime}$ along the direction $\mathbf{n}_{\mathbf{x}}$ as
$\mathbf{x}^{\prime}$ approaches $\mathbf{x}$ and obtain the following
representation involving a hyper-singular kernel,%

\begin{equation}
\frac{\partial}{\partial\mathbf{n}_{\mathbf{x}}}u(\mathbf{x})=-\underset
{\mathbf{x}^{\prime}\rightarrow\mathbf{x}}{\lim}%
{\displaystyle\int\nolimits_{\Gamma\cup S_{a}}}
\frac{\partial^{2}g(\mathbf{x}^{\prime}\mathbf{,y})}{\partial\mathbf{n}%
_{\mathbf{x}}\partial\mathbf{n}_{\mathbf{y}}}u(\mathbf{y})ds,\ \ \mathbf{x}\in
S_{a}.\label{bw39}%
\end{equation}

The integral expression for $\frac{\partial}{\partial\mathbf{n}_{\mathbf{x}}%
}u(\mathbf{x})$ involves two integrals, one regular integral over the upper
hemisphere $\Gamma$ denoted as%

\begin{equation}
\Sigma_{1}=-%
{\displaystyle\int\nolimits_{\Gamma}}
\frac{\partial^{2}g(\mathbf{x,y})}{\partial\mathbf{n}_{\mathbf{x}}%
\partial\mathbf{n}_{\mathbf{y}}}u(\mathbf{y})ds_{y}=-%
{\displaystyle\int\nolimits_{\Gamma}}
\left(  \frac{3}{2\pi}\frac{\cos\theta}{a^{3}}\right)  u(\mathbf{y}%
)ds_{y},\label{bw41}%
\end{equation}
where (\ref{bw33}) has been used in the second equality, and another
hyper-singular integral over the disk $S_{a}$ denoted as%

\begin{equation}
\Sigma_{2}=-\underset{\mathbf{x}^{\prime}\rightarrow\mathbf{x}}{\lim}%
{\displaystyle\int\nolimits_{S_{a}}}
\frac{\partial^{2}g(\mathbf{x}^{\prime}\mathbf{,y})}{\partial\mathbf{n}%
_{\mathbf{x}}\partial\mathbf{n}_{\mathbf{y}}}u(\mathbf{y})ds_{y},\label{bw43}%
\end{equation}
and we have%
\begin{equation}
\frac{\partial}{\partial\mathbf{n}_{\mathbf{x}}}u(\mathbf{x})=\Sigma
_{1}+\Sigma_{2}.\label{bw45}%
\end{equation}

Equation (\ref{bw45}) will be the starting point for the proposed hybrid
method. In computing the integral $\Sigma_{1}$, say by a Gauss quadrature over
the hemisphere surface, we will need the potential solution $u(\mathbf{y})$
for $\mathbf{y}\in\Gamma$ and this solution will be readily computed with the
Feynman-Kac formula (\ref{bw13}) with the WOS as the sampling technique for
the Brownian paths. On the other hand, the singular integral $\Sigma_{2}$,
with appropriate treatment of the hyper-singularities to be described in
detail in the numerical test section, can be calculated directly with the
given Dirichlet boundary data $u(\mathbf{y}),\mathbf{y}\in S_{a}.$ Therefore,
an algorithm using (\ref{bw39}) involves the hybridization of a random walk on
spheres (WOS) and a deterministic boundary integral equation (BIE), which is
termed the \textit{BIE-WOS method}.\bigskip

\noindent\textbf{Remark:} In comparing the last-passage method (\ref{bw29})
with the BIE-WOS method (\ref{bw45}), the former uses the relation between the
Brownian motion of diffusive particles and electric potential from charges on
a conducting surface to arrive at an expression for the surface charge density
based on (\ref{bw29}). On the other hand, the BIE-WOS method uses a
hyper-singular boundary integral equation to get a similar expression in
(\ref{bw45}), which has an additional contribution from the variable potential
on the charged surface (the integral term $\Sigma_{2}$). Both methods use WOS
for particles starting on the hemisphere, but, at different locations. The
last-passage method proposed in \cite{given02} initiates particles' walk
starting from positions all over the hemisphere sampled using a probability
given by (\ref{bw33}), while the BIE-WOS method initiates many particle walks
starting from selected Gauss quadrature points (up to $30\times30$ in our test
problems). Numerical results will show that for problems suitable for both
methods, the total number of particle walk paths and the accuracy and
computational costs are comparable (refer to Test 4 in Section 5.1.3).

\section{Finding Neumann data over a patch of general boundary}

\begin{figure}[ptb]
\begin{center}
\includegraphics[width=3.5in ]{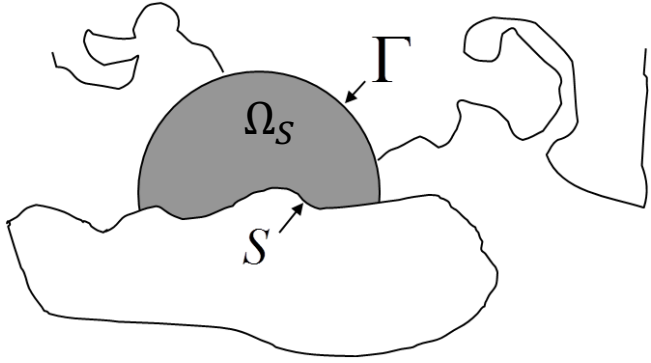}
\end{center}
\caption{Setup of the BIE-WOS method for finding the Neumann data on a patch
$S\subset\partial\Omega.$}%
\label{fig_biewos}%
\end{figure}

In this section, we will extend the BIE-WOS method to the case of general
Dirichlet boundary data and curved domain boundaries. To achieve this goal, we
will superimpose a hemisphere over any selected portion of the boundary
$\partial\Omega$ and denote the intersection portion of the domain boundary by
$S$ and the surface of the hemisphere outside the domain $\Omega$ still by
$\Gamma,$ and the region bounded by $S$ and $\Gamma$ is denoted as $\Omega
_{S}$ (see Fig. \ref{fig_biewos})$.$ Now let $G(\mathbf{x,y})$ be Green's
function of a \textit{whole sphere} with a homogeneous boundary condition,
which can be easily obtained by one Kelvin image charge as discussed before.
Then, the integral representation (\ref{bw03}) can be applied to the boundary
of the domain $\Omega_{S}$ to yield the following identity%

\begin{align}
u(\mathbf{x}) =  & -{\displaystyle\int\nolimits_{\Gamma}} \frac{\partial
G(\mathbf{x,y})}{\partial\mathbf{n}_{\mathbf{y}}}u(\mathbf{y})ds_{y}%
\nonumber\\
+  &
{\displaystyle\int\nolimits_{S}}
\left[  -\frac{\partial G(\mathbf{x,y})}{\partial\mathbf{n}_{\mathbf{y}}%
}u(\mathbf{y})+G(\mathbf{x,y})\frac{\partial u(\mathbf{y})}{\partial
\mathbf{n}_{\mathbf{y}}}\right]  ds_{y},\text{ \ \ }\mathbf{x} \in\Omega
_{S}.\label{bw47}%
\end{align}
It should be noted that the integral over $\Gamma$ only involves the normal
derivative of Green's function as $G$ vanishes on $\Gamma$ by construction. As
a result, only solution $u(\mathbf{y})$ is needed on $\Gamma$ while both
$u(\mathbf{y})$ and the normal derivative $\frac{\partial u(\mathbf{y}%
)}{\partial\mathbf{n}}$ appear in the integral over $S$. As before, the
solution $u(\mathbf{y})$ over $\Gamma$ will be computed with the Feynman-Kac
formula (\ref{bw13}) with WOS and then the Neumann data over $S$ can be solved
from the following integral equation,

\begin{equation}
K\left[  \frac{\partial u}{\partial\mathbf{n}}\right]  (\mathbf{x}%
)=b(\mathbf{x}),\text{ \ \ }\mathbf{x}\in S, \label{bw49}%
\end{equation}
where%
\begin{equation}
K\left[  \frac{\partial u}{\partial\mathbf{n}}\right]  \equiv%
{\displaystyle\int\nolimits_{S}}
G(\mathbf{x,y})\frac{\partial u(\mathbf{y})}{\partial\mathbf{n}_{\mathbf{y}}%
}ds_{y}, \label{bw50}%
\end{equation}
and%
\begin{equation}
b(\mathbf{x})\equiv\left(  \frac{u(\mathbf{x})}{2}+\mathrm{p.v.}%
{\displaystyle\int\nolimits_{S}}
\frac{\partial G(\mathbf{x,y})}{\partial\mathbf{n}_{\mathbf{y}}}%
u(\mathbf{y})ds_{y})\right)  +%
{\displaystyle\int\nolimits_{\Gamma}}
\frac{\partial G(\mathbf{x,y})}{\partial\mathbf{n}_{\mathbf{y}}}%
u(\mathbf{y})ds_{y}, \label{bw51}%
\end{equation}
where p.v. stands for the Cauchy principal value of the double layer potential
\cite{cai13}.

The integral equation (\ref{bw49}) is of the first kind which is
ill-conditioned and may cause numerical difficulties especially when the
algebraic system from discretization becomes large. When that happens, a
well-conditioned second kind of integral equation can be obtained by taking
normal derivative of (\ref{bw47}), resulting in the following identity %

\begin{align}
\frac{\partial}{\partial\mathbf{n}_{\mathbf{x}}}u(\mathbf{x}) &  =-%
{\displaystyle\int\nolimits_{\Gamma}}
\frac{\partial^{2}G(\mathbf{x,y})}{\partial\mathbf{n}_{\mathbf{x}}%
\partial\mathbf{n}_{\mathbf{y}}}u(\mathbf{y})ds_{y}\nonumber\\
&  +%
{\displaystyle\int\nolimits_{S}}
\left[  -\frac{\partial^{2}G(\mathbf{x,y})}{\partial\mathbf{n}_{\mathbf{x}%
}\partial\mathbf{n}_{\mathbf{y}}}u(\mathbf{y})+\frac{\partial G(\mathbf{x,y}%
)}{\partial\mathbf{n}_{\mathbf{x}}}\frac{\partial u(\mathbf{y})}%
{\partial\mathbf{n}_{\mathbf{y}}}\right]  ds_{y},\text{ \ \ }\mathbf{x}%
\in\Omega_{S}.\label{bw53}%
\end{align}
Let $\mathbf{x}$ approach the boundary $S$. We obtain the following second
kind integral equation%

\begin{equation}
\left(  \frac{1}{2}I-D\right)  \left[  \frac{\partial u}{\partial\mathbf{n}%
}\right]  (\mathbf{x})=b(\mathbf{x}),\text{ \ \ }\mathbf{x}\in S, \label{bw55}%
\end{equation}
where the integral operator of a double layer potential is%

\begin{equation}
D\left[  \frac{\partial u}{\partial\mathbf{n}}\right]  (\mathbf{x})\equiv%
{\displaystyle\int\nolimits_{S}}
\frac{\partial G(\mathbf{x,y})}{\partial\mathbf{n}_{\mathbf{x}}}\frac{\partial
u(\mathbf{y})}{\partial\mathbf{n}_{\mathbf{y}}}ds_{y}, \label{bw57}%
\end{equation}
and%

\begin{equation}
b(\mathbf{x})\equiv-%
{\displaystyle\int\nolimits_{\Gamma}}
\frac{\partial^{2}G(\mathbf{x,y})}{\partial\mathbf{n}_{\mathbf{x}}%
\partial\mathbf{n}_{\mathbf{y}}}u(\mathbf{y})ds_{y}-\mathrm{p.f.}%
{\displaystyle\int\nolimits_{S}}
\frac{\partial^{2}G(\mathbf{x,y})}{\partial\mathbf{n}_{\mathbf{x}}%
\partial\mathbf{n}_{\mathbf{y}}}u(\mathbf{y})ds_{y},\text{ \ \ }\mathbf{x}\in
S,\label{bw59}%
\end{equation}
and p.f. denotes the Hadamard finite part limit for the hyper-singular
integral \cite{cai13}, which can be handled by a regularization
technique.\bigskip

\textbf{BIE-WOS Algorithm:} The BIE-WOS method for the Neumann data over a
patch $S$ will consist of two steps:\medskip\bigskip

\begin{itemize}
\item Step 1: Apply the Feynman-Kac formula (\ref{bw13}) with the WOS sampling
technique to compute the potential solution $u(\mathbf{y}_{i,j})$\ at Gauss
points $\mathbf{y}_{i,j}\in\Gamma$. Compute the right-hand-side function
$b(\mathbf{x})$ in (\ref{bw51}) or (\ref{bw59}) by some Gauss
quadratures.\bigskip

\item Step 2: Solve the BIE (\ref{bw49}) or (\ref{bw55}) with a collocation
method for the Neumann data $\frac{\partial u}{\partial\mathbf{n}}$ over $S$.
\end{itemize}

\bigskip

\noindent\textbf{Remark:} To find the derivatives of the potential inside and
outside a bounded domain $\Omega,$ the BIE-WOS method is applied separately.
Namely, the hemisphere $\Omega_{s}$ \ in Fig. \ref{fig_biewos} will be located in the
interior and exterior of $\Omega$ (wherein the WOS method will be used),
respectively. Once the Neumann data for the respective potential is found by
the BIE-WOS method, the integral representation (\ref{bw03}) can be used with
the corresponding Green's function to obtain the potential in the whole
interior and exterior space.

\section{Numerical Results}

In this section, we will present a series of numerical tests to demonstrate
the accuracy and efficiency of the proposed BIE-WOS method for finding the
Neumann data at a single point on a flat boundary or on a patch over a curved boundary.

\subsection{Finding Neumann data at one point on a flat boundary}

\subsubsection{Regularization of hyper-singular integrals}

First, let us present a regularization method using simple solution of the
Laplace equation \cite{nedelec78} to compute the hyper-singular integral in
(\ref{bw43}) and (\ref{bw59}). First, with some simple calculations, the term
$\Sigma_{2}$ of (\ref{bw43}) is found to be a Hadamard finite part limit of
the following hyper-singular integral:%
\begin{equation}
\Sigma_{2}=-\underset{\mathbf{x}^{\prime}\rightarrow\mathbf{x}}{\lim}%
{\displaystyle\int\nolimits_{S_{a}}}
\frac{\partial^{2}g(\mathbf{x}^{\prime}\mathbf{,y})}{\partial\mathbf{n}%
_{\mathbf{x}}\partial\mathbf{n}_{\mathbf{y}}}u(\mathbf{y})ds_{y}%
=-\mathrm{p.f.}%
{\displaystyle\int\nolimits_{S_{a}}}
\frac{1}{2\pi}\left(  \frac{1}{\rho^{3}}-\frac{1}{a^{3}}\right)
u(\mathbf{y})ds_{y},\label{bw61}%
\end{equation}
where $\rho=|\mathbf{x-y}|,\mathbf{x},\mathbf{y\in}$ $S_{a}.$ The finite part
(p.f.) limit of Hadamard type is defined by removing a divergent part in the
process of defining a principal value (i.e. by removing a small patch of size
$\varepsilon$ centered at $\mathbf{x}$ and then let $\varepsilon$ approach
zero) \cite{cai13}. For the Laplace equation considered here, we can
regularize this hyper-singularity by invoking an integral identity for the
special solution $u\equiv\phi(\mathbf{x})$, with $\mathbf{x}$ being fixed,
namely, the integral identity (\ref{bw39}) applied to this constant solution
results in%

\begin{equation}
0=-%
{\displaystyle\int\nolimits_{\Gamma}}
\frac{\partial^{2}g(\mathbf{x,y})}{\partial\mathbf{n}_{\mathbf{x}}%
\partial\mathbf{n}_{\mathbf{y}}}\phi(\mathbf{x})ds-\underset{\mathbf{x}%
^{\prime}\rightarrow\mathbf{x}}{\lim}%
{\displaystyle\int\nolimits_{S_{a}}}
\frac{\partial^{2}g(\mathbf{x}^{\prime}\mathbf{,y})}{\partial\mathbf{n}%
_{\mathbf{x}}\partial\mathbf{n}_{\mathbf{y}}}\phi(\mathbf{x})ds,\text{
\ }\mathbf{x}\in S.\label{bw63}%
\end{equation}

Subtracting (\ref{bw63}) from (\ref{bw45}), we have a modified formula for the
Neumann data as%

\begin{equation}
\frac{\partial}{\partial\mathbf{n}_{\mathbf{x}}}u(\mathbf{x})=\Sigma
_{1}^{\prime}+\Sigma_{2}^{\prime},\text{ \ \ \ }\mathbf{x}\in S, \label{bw64}%
\end{equation}
where $\Sigma_{1}^{\prime}$ and $\Sigma_{2}^{\prime}$ are now regularized
versions of $\Sigma_{1}$ and $\Sigma_{2}$\ in (\ref{bw41}) and (\ref{bw43}),
respectively, i.e.,%

\begin{equation}
\Sigma_{1}^{\prime}=-%
{\displaystyle\int\nolimits_{\Gamma}}
\frac{\partial^{2}g(\mathbf{x,y})}{\partial\mathbf{n}_{\mathbf{x}}%
\partial\mathbf{n}_{\mathbf{y}}}\left(  u(\mathbf{y})-\phi(\mathbf{x})\right)
ds_{y},\label{bw65}%
\end{equation}
and%
\begin{align}
\Sigma_{2}^{\prime} &  =-\underset{\mathbf{x}^{\prime}\rightarrow\mathbf{x}%
}{\lim}%
{\displaystyle\int\nolimits_{S_{a}}}
\frac{\partial^{2}g(\mathbf{x}^{\prime}\mathbf{,y})}{\partial\mathbf{n}%
_{\mathbf{x}}\partial\mathbf{n}_{\mathbf{y}}}\left(  u(\mathbf{y}%
)-\phi(\mathbf{x})\right)  ds_{y}\nonumber\\
&  =-\underset{\mathbf{x}^{\prime}\rightarrow\mathbf{x}}{\lim}%
{\displaystyle\int\nolimits_{S_{a}}}
\frac{1}{2\pi}\left(  \frac{1}{r^{3}}-\frac{1}{a^{3}}\right)  \left(
\phi(\mathbf{y})-\phi(\mathbf{x})\right)  ds_{y},\label{bw68}%
\end{align}
where $\mathbf{x}^{\prime}=\mathbf{x+}(0,0,\varepsilon),$ $r=\sqrt{\rho
^{2}+\varepsilon^{2}},\boldsymbol{\rho }=\mathbf{x-y,}\rho=|\mathbf{x-y}%
|,\mathbf{x,y}\in S_{a}$. Moreover, the boundary condition $u(\mathbf{y}%
)=\phi(\mathbf{y})$, $\mathbf{y}\in S_{a}$ has been invoked in (\ref{bw68}).

Compared with (\ref{bw61}), the singularity in the integral $\Sigma
_{2}^{\prime}$ in (\ref{bw68}) has been weakened by the factor $\left(
\phi(\mathbf{y})-\phi(\mathbf{x})\right)  ,$ which vanishes at $\mathbf{x,}$
and $\Sigma_{2}^{\prime}$ will be evaluated by a Gauss quadrature. Let us only
consider the integral involving the singular term $\frac{1}{r^{3}}$ in
(\ref{bw68}), which is denoted by $\Sigma_{2}^{\ast},$ i.e.,%

\begin{equation}
\Sigma_{2}^{\ast}=-\frac{1}{2\pi}\underset{\mathbf{x}^{\prime}\rightarrow
\mathbf{x}}{\lim}%
{\displaystyle\int\nolimits_{S_{a}}}
\frac{1}{r^{3}}\left(  \phi(\mathbf{y})-\phi(\mathbf{x})\right)  ds_{y}.
\label{bw71}%
\end{equation}

Consider a circular patch $\Lambda_{\delta}$ of radius $\delta$ centered at
$\mathbf{x}$, and then $\Sigma_{2}^{\ast}$ can be split further into two
integrals as follows%
\begin{align}
\Sigma_{2}^{\ast}  &  =-\frac{1}{2\pi}%
{\displaystyle\int\nolimits_{S_{a}\backslash\Lambda_{\delta}}}
\frac{1}{\rho^{3}}\left(  \phi(\mathbf{y})-\phi(\mathbf{x})\right)
ds_{y}\nonumber\\
&  \text{ \ \ }-\frac{1}{2\pi}\underset{\mathbf{x}^{\prime}\rightarrow
\mathbf{x}}{\lim}%
{\displaystyle\int\nolimits_{\Lambda_{\delta}}}
\frac{1}{r^{3}}\left(  \phi(\mathbf{y})-\phi(\mathbf{x})\right)
ds_{y}\nonumber\\
&  =-\frac{1}{2\pi}%
{\displaystyle\int\nolimits_{S_{a}\backslash\Lambda_{\delta}}}
\frac{1}{\rho^{3}}\left(  \phi(\mathbf{y})-\phi(\mathbf{x})\right)
ds_{y}+\Delta. \label{bw73}%
\end{align}

To estimate the term $\Delta$, we apply a Taylor expansion of the boundary
data $\phi(\mathbf{y})$ at $\mathbf{x}$%

\begin{equation}
\phi(\mathbf{y})-\phi(\mathbf{x})=\nabla\phi(\mathbf{x})\cdot
\boldsymbol{\rho }+O(\rho^{2}). \label{bw75}%
\end{equation}
Then, we obtain%
\begin{align}
\Delta &  =-\frac{\nabla\phi(\mathbf{x})\cdot}{2\pi}\underset{\mathbf{x}%
^{\prime}\rightarrow\mathbf{x}}{\lim}%
{\displaystyle\int\nolimits_{\Lambda_{\delta}}}
\frac{\boldsymbol{\rho}}{r^{3}}ds_{y}+\frac{1}{2\pi}%
{\displaystyle\int\nolimits_{\Lambda_{\delta}}}
\frac{O(\rho^{2})}{r^{3}}ds_{y}\nonumber\\
&  =-\frac{\nabla\phi(\mathbf{x})\cdot}{2\pi}\underset{\mathbf{x}^{\prime
}\rightarrow\mathbf{x}}{\lim}%
{\displaystyle\int\nolimits_{0}^{\delta}}
{\displaystyle\int\nolimits_{0}^{2\pi}}
\frac{\rho\mathbf{(}\cos\theta,\sin\theta)}{(\rho^{2}+\varepsilon^{2})^{3/2}%
}\rho d\theta d\rho\nonumber\\
&  \text{ \ \ }+\text{ }\underset{\mathbf{x}^{\prime}\rightarrow\mathbf{x}%
}{\lim}%
{\displaystyle\int\nolimits_{0}^{\delta}}
\frac{O(\rho^{2})}{(\rho^{2}+\varepsilon^{2})^{3/2}}\rho d\rho\nonumber\\
&  =0+\underset{\mathbf{x}^{\prime}\rightarrow\mathbf{x}}{\lim}%
{\displaystyle\int\nolimits_{0}^{\delta}}
\frac{O(\rho^{3})}{(\rho^{2}+\varepsilon^{2})^{3/2}}d\rho. \label{bw77}%
\end{align}
Now for all positive $\varepsilon>0,$ we have
\begin{equation}
\frac{\rho^{3}}{(\rho^{2}+\varepsilon^{2})^{3/2}}\leq1. \label{bw79}%
\end{equation}
As a result, the following estimate of the term $\Delta$ holds%

\begin{equation}
\Delta=O(\delta). \label{bw81}%
\end{equation}
Finally, the regularized integral $\Sigma_{2}^{\ast}$ will be approximated by
the integral over $S_{a}\backslash\Lambda_{\delta}$ with an accuracy of
$O(\delta)$ and a Gauss quadrature formula over the ring shaped region
$S_{a}\backslash\Lambda_{\delta}$:%

\begin{equation}
\Sigma_{2}^{\ast}=-\frac{1}{2\pi}%
{\displaystyle\int\nolimits_{S_{a}\backslash\Lambda_{\delta}}}
\frac{1}{\rho^{3}}\left(  \phi(\mathbf{y})-\phi(\mathbf{x})\right)
ds_{y}+O(\delta). \label{bw82}%
\end{equation}

\subsubsection{Gauss quadratures over the hemisphere $\Gamma$ and
$S_{a}\backslash\Lambda_{\delta}$ and WOS}

To compute the integral $\Sigma_{1}^{\prime}$, we use $N_{g1}\times N_{g1}$
Gauss points over the hemispherical surface $\Gamma$%
\begin{equation}
\Sigma_{1}^{\prime}\simeq-%
{\displaystyle\sum\limits_{i,j=1}^{N_{g1}}}
\omega_{i}\omega_{j}\frac{\pi^{2}}{4}(a^{2}\sin\theta_{i})\frac{3}{2a}\left(
\frac{\cos\theta_{i}}{\pi a^{2}}\right)  \left(  u(\mathbf{y}_{i,j}%
)-\phi(\mathbf{x})\right)  , \label{bw83}%
\end{equation}
where%
\begin{equation}
\theta_{i}=\frac{\pi}{4}(\xi_{i}+1),\varphi_{j}=\pi(\xi_{j}+1),\mathbf{y}%
_{i,j}=(a,\theta_{i},\varphi_{j}), \label{bw85}%
\end{equation}
and $\omega_{i}$ and $\xi_{i},1\leq i\leq$ $N_{g1}$ are the Gauss quadrature
weights and locations, respectively. $\frac{\pi^{2}}{4}(a^{2}\sin\theta_{i})$
is the area of the surface element in the spherical coordinates.

Now, each of the solution values $u(\mathbf{y}_{i,j}),$ $\mathbf{y}_{i,j}%
\in\Gamma$ will be obtained by the Feynman-Kac formula (\ref{bw13}) with
$N_{path}$ Brownian particles all starting from $\mathbf{y}_{i,j}$, namely%
\begin{equation}
u(\mathbf{y}_{i,j})\simeq\frac{1}{N_{path}}%
{\displaystyle\sum\limits_{\text{ k=1}}^{N_{path}}}
\phi(\mathbf{e}_{k}), \label{bw87}%
\end{equation}
where $\mathbf{e}_{k}$ is the location on $\partial\Omega$ where a path terminates.

The total number $N_{\mathrm{path-bie-wos}}$ of Brownian particles needed in
the BIE-WOS method will be%
\begin{equation}
N_{\mathrm{path-bie-wos}}=N_{g1}\times N_{g1}\times N_{path}. \label{bw88}%
\end{equation}

Next, the integral $\Sigma_{2}^{\ast}$ in (\ref{bw82}) will be computed with
another $N_{g2}\times N_{g2}$ Gauss quadrature over the ring shaped region
$S_{a}\backslash\Lambda_{\delta}$ with an error of $O(\delta)$ in addition to
the error from the Gauss quadrature.

\subsubsection{Numerical tests}

In this section, we will present several numerical tests to demonstrate the
accuracy and efficiency of the proposed BIE-WOS method for finding the Neumann
data at a given point over a flat boundary for general Dirichlet boundary
data. For comparison, we also implement the last-passage Monte Carlo method
proposed in \cite{given02}. For accuracy comparison, the charge density is
calculated with the FastCap, an open-source code developed in MIT
\cite{white91} for 3-D capacitance extraction tool in industry and academia.
The Fastcap is an indirect BEM, accelerated by the fast multipole method
(FMM), and its linear system is solved by a conjugate gradient method. For the
case of complex potentials on the surfaces, we also implemented a direct BEM
(DBEM) \cite{yu04}. To identify Brownian particles going to infinity, a large
sphere of radius of $10^{5}$ is used, which is found to be large enough for
the desired accuracy. Thus, once a particle gets out this sphere, it will be
considered as having gone to infinity.

\bigskip

\noindent$\bullet$ \textit{Test 1- Charge densities on a planar interface
between two dielectric half spaces} \bigskip

\begin{figure}[ptb]
\begin{center}
\includegraphics[width=3.5in ]{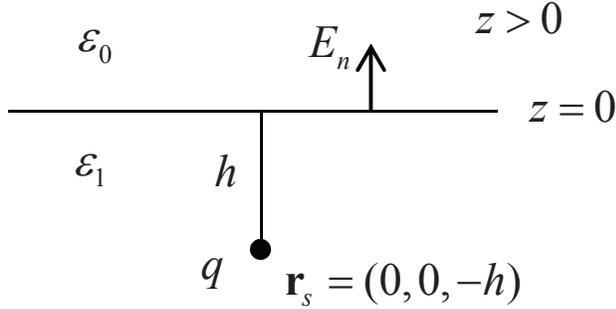}
\end{center}
\caption{Potential above a half-space}%
\label{fig_halfspace}%
\end{figure}

As shown in Fig. \ref{fig_halfspace}, the whole space is divided by a planar
interface between two dielectric domains, and the dielectric constants are
$\epsilon_{0}$ and $\epsilon_{1}$ in the upper and lower domain, respectively.
A charge $q$ is located at $\mathbf{r}_{s}=(0,0,-h)$. Then the potential in
the upper space is given by%

\begin{equation}
u(\mathbf{r})=\frac{q^{\prime}}{4\pi\epsilon_{0}}\frac{1}{|\mathbf{r}%
-\mathbf{r}_{s}|},\ \ q^{\prime}=\frac{2\epsilon_{0}}{\epsilon_{0}%
+\epsilon_{1}}q, \label{bw89}%
\end{equation}
and $u(\mathbf{r})$ satisfies the Laplace equation $\nabla^{2}u(\mathbf{r}%
)=0,z>0$ with a variable Dirichlet data on the boundary $z=0$.

\bigskip

The charge density at the point $\mathbf{x}=(0.5,0,0)$ by the last-passage
method and the BIE-WOS method with various radius $a$ of the hemisphere are
listed in Table \ref{tab_halfspace}. In the last-passage method, the total
number of the sampling paths is $N=4\times10^{5}$. In the BIE-WOS method, the
number of Gauss points is $N_{g1}\times N_{g1}=20\times20$ for the hemisphere,
and is $N_{g2}\times N_{g2}=20\times20$ for the integral on the 2-D disk
$S_{a}$. Starting from each Gauss point on the hemisphere, the number of the
sampling paths is $N_{path}=10^{3}$. Therefore, the total number of paths for
the BIE-WOS method is also $4\times10^{5}$. In both methods, the thickness
$\varepsilon$ of the absorption layer for the WOS method is taken to be
$10^{-5}$.

\begin{table}[ptb]
\caption{Charge density on the planar interface with different radius}%
\label{tab_halfspace}
\begin{center}%
\begin{tabular}
[c]{|l|c|c|c|c|c|c|l|}\hline
& \multicolumn{2}{|l|}{\ \ Last-passage} &
\multicolumn{4}{|l|}{\ \ \ \ \ \ \ \ \ \ \ \ \ BIE-WOS} &
analytical\\\cline{2-7}%
$a$ & $\Sigma_{\text{LP}}$ & err$\%$ & $\Sigma_{1}^{^{\prime}}$ &
$\ \ \Sigma_{2}^{^{\prime}}$ & $\Sigma_{1}^{^{\prime}}+\Sigma_{2}^{^{\prime}}$
& err$\%$ & \multicolumn{1}{|c|}{solution}\\\hline
0.1 & \multicolumn{1}{|l|}{0.698543} & \multicolumn{1}{|r|}{-2.38} &
\multicolumn{1}{|l|}{0.69884} & \multicolumn{1}{|l|}{0.018777} &
\multicolumn{1}{|l|}{0.717612} & \multicolumn{1}{|r|}{0.29} & \\\cline{1-7}%
0.2 & \multicolumn{1}{|l|}{0.677996} & \multicolumn{1}{|r|}{-5.25} &
\multicolumn{1}{|l|}{0.67784} & \multicolumn{1}{|l|}{0.037515} &
\multicolumn{1}{|l|}{0.715355} & \multicolumn{1}{|r|}{-0.03} & \\\cline{1-7}%
0.5 & \multicolumn{1}{|l|}{0.622949} & \multicolumn{1}{|r|}{-12.94} &
\multicolumn{1}{|l|}{0.62146} & \multicolumn{1}{|l|}{0.093054} &
\multicolumn{1}{|l|}{0.714517} & \multicolumn{1}{|r|}{-0.14} &
0.71554\\\cline{1-7}%
0.7 & \multicolumn{1}{|l|}{0.586721} & \multicolumn{1}{|r|}{-18.00} &
\multicolumn{1}{|l|}{0.58432} & \multicolumn{1}{|l|}{0.128971} &
\multicolumn{1}{|l|}{0.713287} & \multicolumn{1}{|r|}{-0.32} & \\\cline{1-7}%
1.0 & \multicolumn{1}{|l|}{0.534695} & \multicolumn{1}{|r|}{-25.27} &
\multicolumn{1}{|l|}{0.53659} & \multicolumn{1}{|l|}{0.179973} &
\multicolumn{1}{|l|}{0.716562} & \multicolumn{1}{|r|}{0.14} & \\\hline
\end{tabular}
\end{center}
\end{table}

From Table \ref{tab_halfspace}, we can see that when the radius increases, the
relative error of the last-passage method grows and grows even up to
$-25.27\%$. It shows that when the potential Dirichlet data on the disk
$S_{a}$ is not constant, the last-passage method is not applicable. The
variable potential inside the disk $S_{a}$ will influence the charge density
at $\mathbf{x}$. In contrast, the BIE-WOS method includes such influences as
shown in the results, and most importantly, is independent of the radius $a$,
for its maximal relative error is less than 0.32\% when the radius ranges from
0.1 to 1.0.

\begin{table}[ptb]
\caption{Accuracy of the de-singularization in (4.12)}%
\label{tab_accgauss}
\begin{center}%
\begin{tabular}
[c]{|l|l|l|l|l|l|l|l|l|}\hline
$\delta/a$ & \multicolumn{8}{|l|}{$\ \ \ \ \ \ \Sigma_{2}^{\prime}%
\ \ $calculated by$\ N_{g2}\times N_{g2}$ Gauss Quadrature}\\\cline{2-9}
& $4\times4$ & err$\%$ & $6\times6$ & err$\%$ & $10\times10$ & err$\%$ &
$20\times20$ & err$\%$\\\hline
\multicolumn{1}{|r|}{$10^{-1}$} & \multicolumn{1}{|r|}{0.09659} &
\multicolumn{1}{|r|}{3.800} & \multicolumn{1}{|r|}{0.08983} &
\multicolumn{1}{|r|}{-3.462} & \multicolumn{1}{|r|}{0.08949} &
\multicolumn{1}{|r|}{-3.833} & \multicolumn{1}{|r|}{0.08949} &
\multicolumn{1}{|r|}{-3.8351}\\\hline
\multicolumn{1}{|r|}{$10^{-2}$} & \multicolumn{1}{|r|}{0.10042} &
\multicolumn{1}{|r|}{7.916} & \multicolumn{1}{|r|}{0.09306} &
\multicolumn{1}{|r|}{0.001} & \multicolumn{1}{|r|}{0.09273} &
\multicolumn{1}{|r|}{-0.347} & \multicolumn{1}{|r|}{0.09273} &
\multicolumn{1}{|r|}{-0.3492}\\\hline
\multicolumn{1}{|r|}{$10^{-3}$} & \multicolumn{1}{|r|}{0.10083} &
\multicolumn{1}{|r|}{8.352} & \multicolumn{1}{|r|}{0.09335} &
\multicolumn{1}{|r|}{0.314} & \multicolumn{1}{|r|}{0.09302} &
\multicolumn{1}{|r|}{-0.033} & \multicolumn{1}{|r|}{0.09302} &
\multicolumn{1}{|r|}{-0.0346}\\\hline
\multicolumn{1}{|r|}{$10^{-4}$} & \multicolumn{1}{|r|}{0.10087} &
\multicolumn{1}{|r|}{8.397} & \multicolumn{1}{|r|}{0.09337} &
\multicolumn{1}{|r|}{0.345} & \multicolumn{1}{|r|}{0.09305} &
\multicolumn{1}{|r|}{-0.002} & \multicolumn{1}{|r|}{0.09305} &
\multicolumn{1}{|r|}{-0.0035}\\\hline
\multicolumn{1}{|r|}{$10^{-5}$} & \multicolumn{1}{|r|}{0.10087} &
\multicolumn{1}{|r|}{8.402} & \multicolumn{1}{|r|}{0.09338} &
\multicolumn{1}{|r|}{0.348} & \multicolumn{1}{|r|}{0.09306} &
\multicolumn{1}{|r|}{0.001} & \multicolumn{1}{|r|}{0.09305} &
\multicolumn{1}{|r|}{-0.0003}\\\hline
\multicolumn{1}{|r|}{$10^{-6}$} & \multicolumn{1}{|r|}{0.10087} &
\multicolumn{1}{|r|}{8.402} & \multicolumn{1}{|r|}{0.09338} &
\multicolumn{1}{|r|}{0.348} & \multicolumn{1}{|r|}{0.09306} &
\multicolumn{1}{|r|}{0.001} & \multicolumn{1}{|r|}{0.09305} &
\multicolumn{1}{|r|}{}\\\hline
\end{tabular}
\end{center}
\end{table}

Table \ref{tab_accgauss} lists the accuracy of the de-singularized $\Sigma
_{2}^{\prime}$ in (\ref{bw82}) with different values $\delta$ and numbers of
Gauss points $N_{g2}\times N_{g2}$, where the location of the sought-after
density is at $(0.5,0,0)$. The result of $N_{g2}\times N_{g2}=20\times20$ with
$\delta/a=10^{-6}$ is taken as the reference value for $\Sigma_{2}^{\prime}$.
Table \ref{tab_accgauss} shows the convergence speed of $\Sigma_{2}^{\prime}$
as $\delta/a$ goes to zero and the number of the Gauss points increases. It
can be seen that when the number of the Gauss points is large enough, for
example $20\times20$, the relative error is on the order of $\delta/a$,
verifying the estimate in (\ref{bw82}) .

\bigskip

\noindent$\bullet$ \textit{Test 2: Four rectangular plates with a piecewise
constant potential distribution}

\bigskip

\begin{figure}[ptb]
\begin{center}
\includegraphics[width=3.5in ]{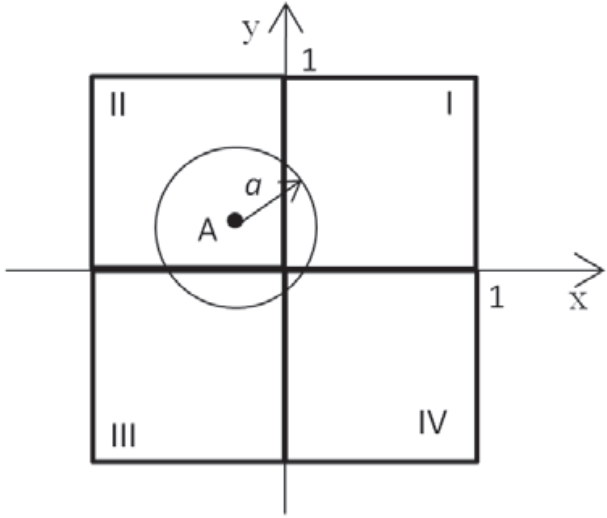}
\end{center}
\caption{Four plates at different potentials}%
\label{fig_fourplates}%
\end{figure}

A 3-D structure with four rectangular plates is depicted in Fig.
\ref{fig_fourplates}, where the length, width and thickness of all four plates
are $1\mathrm{m}\times1\mathrm{m}\times0.01\mathrm{m}$. \ First, we set the
potential of plate II to 1V and the potential of the other three (I, III and
IV) to 0V, and compute the charge density at the point $A(-0.2273,0.2273)$.
The results of all four methods are listed in Table \ref{tab_fourplates_01},
taking the results by the FastCap as the reference where each side of the
plates is discretized into $99\times99$ panels. The DBEM uses a discretization
with $11\times11$ panels on each side, and its relative error is -0.46\%.

Both the last-passage method and the BIE-WOS method run with
various radius $a$ of the hemisphere, and the parameters are the same as in
Test 1. In this case, the integral $\Sigma_{2}^{\prime}$ is related to the
area of the intersecting area between the disk $S_{a}$ \ and the plates I, III
and IV, and we just compute it directly by the \textit{quad} function in
Matlab, instead of Gauss quadratures.

Note that the potential on the boundary $\partial\Omega$ here is piecewise
constant. Therefore, in the last-passage method, charge density should be
computed, instead of by (\ref{bw35}), by the following formula:\qquad\qquad%
\begin{equation}
\Sigma_{\mathrm{LP}}=\frac{3}{2a}\frac{N_{\inf}+N_{I}+N_{III}+N_{IV}%
}{N_{\mathrm{path-LP}}}, \label{bw91}%
\end{equation}
where $N_{\inf},N_{I},N_{III},$ and $N_{IV}$ represent the number of particles
which finally go to infinity, plate I, III, and IV, respectively.
$N_{\mathrm{path-LP}}$ denotes the total number of Brownian paths starting
from the hemisphere $\Gamma.$

\begin{table}[ptb]
\caption{Charge density of a structure of four unit plates with different
radius}%
\label{tab_fourplates_01}
\begin{center}%
\begin{tabular}
[c]{|l|c|c|c|c|c|c|c|c|}\hline
$a$ & \multicolumn{2}{|l}{Last-passage} &
\multicolumn{4}{|l}{\ \ \ \ \ \ \ \ BIE-WOS} & \multicolumn{2}{|l|}{DBEM}%
\\\cline{2-9}
& $\Sigma_{\text{LP}}$ & err$\%$ & $\Sigma_{1}^{^{\prime}}$ & $\Sigma
_{2}^{^{\prime}}$ & $\Sigma_{1}^{^{\prime}}+\Sigma_{2}^{^{\prime}}$ & err$\%$
& value & err$\%$\\\hline
\multicolumn{1}{|r|}{0.1} & \multicolumn{1}{|r|}{2.6084} &
\multicolumn{1}{|r|}{0.05} & \multicolumn{1}{|r|}{2.6051} &
\multicolumn{1}{|r|}{0} & \multicolumn{1}{|r|}{2.6051} &
\multicolumn{1}{|r|}{-0.07} & \multicolumn{1}{|r|}{} & \multicolumn{1}{|r|}{}%
\\\cline{1-7}%
\multicolumn{1}{|r|}{0.2} & \multicolumn{1}{|r|}{2.6026} &
\multicolumn{1}{|r|}{-0.17} & \multicolumn{1}{|r|}{2.6051} &
\multicolumn{1}{|r|}{0} & \multicolumn{1}{|r|}{2.6051} &
\multicolumn{1}{|r|}{-0.07} & \multicolumn{1}{|r|}{2.595} &
\multicolumn{1}{|r|}{-0.46}\\\cline{1-7}%
\multicolumn{1}{|r|}{0.2273} & \multicolumn{1}{|r|}{2.6099} &
\multicolumn{1}{|r|}{0.11} & \multicolumn{1}{|r|}{2.6064} &
\multicolumn{1}{|r|}{0} & \multicolumn{1}{|r|}{2.6064} &
\multicolumn{1}{|r|}{-0.02} & \multicolumn{1}{|r|}{} & \multicolumn{1}{|r|}{}%
\\\hline
\multicolumn{1}{|r|}{0.3} & \multicolumn{1}{|r|}{2.5252} &
\multicolumn{1}{|r|}{-3.14} & \multicolumn{1}{|r|}{2.5178} &
\multicolumn{1}{|r|}{0.0892} & \multicolumn{1}{|r|}{2.6070} &
\multicolumn{1}{|r|}{-0.00} & \multicolumn{2}{|c|}{Fastcap}\\\hline
\multicolumn{1}{|r|}{0.5} & \multicolumn{1}{|r|}{1.9698} &
\multicolumn{1}{|r|}{-24.44} & \multicolumn{1}{|r|}{1.9692} &
\multicolumn{1}{|r|}{0.6330} & \multicolumn{1}{|r|}{2.6022} &
\multicolumn{1}{|r|}{-0.19} & \multicolumn{2}{|c|}{}\\\cline{1-7}%
\multicolumn{1}{|r|}{0.7} & \multicolumn{1}{|r|}{1.5779} &
\multicolumn{1}{|r|}{-39.48} & \multicolumn{1}{|r|}{1.5784} &
\multicolumn{1}{|r|}{1.0271} & \multicolumn{1}{|r|}{2.6055} &
\multicolumn{1}{|r|}{-0.06} & \multicolumn{2}{|c|}{2.607}\\\hline
\end{tabular}
\end{center}
\end{table}

From Table \ref{tab_fourplates_01}, we can see that when the radius
$a\leq0.2773$, i.e. the disk $S_{a}$ is totally inside the plate II, the
last-passage method is correct with a maximal relative error less than 0.17\%.
However, once $S_{a}$ becomes larger and covers areas of plates with different
potentials, the relative error of the last-passage method increases, even up
to $-39.48\%$. In comparison, the BIE-WOS method maintains its accuracy
insensitive to the radius $a$ with a maximal relative errors less than
$-0.19\%$ as the radius varies from 0.1 to 0.7. This again confirms the fact
that the last-passage method of \cite{given02} is designed for conducting
surfaces (i.e., constant potential), not for surfaces of variable potentials.
Therefore, it should not be used when the disk $S_{a}$ includes regions of
different potentials.

In conclusion, for a general variable potential, the last-passage method is
limited while the BIE-WOS method does not suffer from the constraint of a
constant boundary potential.

\bigskip

\noindent$\bullet$ \textit{Test 3: Four rectangular plates with a complex
potential distribution }

\bigskip

To further emphasize the point raised above in Test 2, we set the four plates
with a complex potential distribution as:
\begin{equation}
\phi(x,y)=\sin mx\sin ny. \label{bw93}%
\end{equation}
To obtain an accurate result, the last-passage method will require
increasingly smaller radius $a$ for ever larger $m$ and $n$ to achieve an
(approximately) constant potential within the disk $S_{a}.$

The charge density at the point $(-0.5,0.5)$ by the last-passage method, the
BIE-WOS method and the DBEM are shown in Table \ref{tab_fourplates_sin}. We
take the result of the DBEM with $17\times17$ panels on each plate as the
reference solution. All other parameters in the BIE-WOS and last-passage
methods are same as in the previous case. From Table \ref{tab_fourplates_sin},
we can see that the BIE-WOS method is more accurate.

\begin{figure}[ptb]
\begin{center}
\includegraphics[width=4in ]{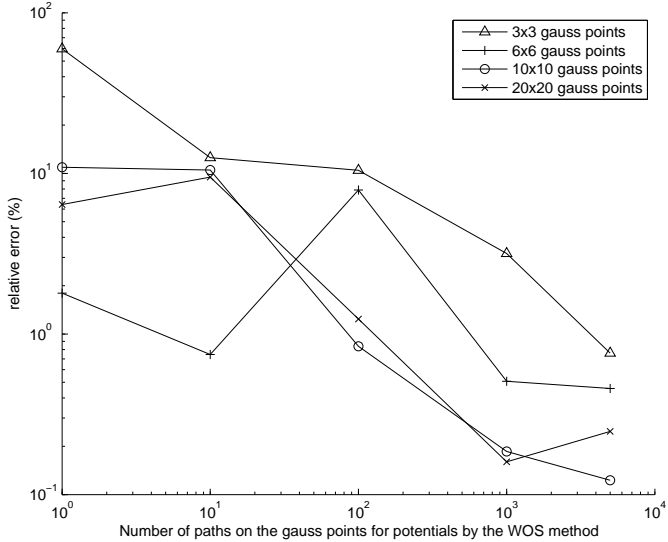}
\end{center}
\caption{Convergence of the BIE-WOS method vs the number of Brownian paths and
Gauss points}%
\label{fig_pathcomp}%
\end{figure}

\begin{table}[ptb]
\caption{Charge density of a structure of four unit plates with complex
voltages in different radius}
\begin{center}%
\begin{tabular}
[c]{|r|c|c|c|c|c|c|c|}\hline
$a$ & \multicolumn{2}{|c|}{Last--passage} & \multicolumn{4}{|c|}{BIE-WOS} &
DBEM\\\cline{2-8}
& $\Sigma_{\text{LP}}$ & err$\%$ & $\Sigma_{1}^{^{\prime}}$ & $\Sigma
_{2}^{^{\prime}}$ & $\Sigma_{1}^{^{\prime}}+\Sigma_{2}^{^{\prime}}$ & err$\%$
& value\\\hline
0.1 & \multicolumn{1}{|r|}{-0.4522} & \multicolumn{1}{|r|}{3.92} &
\multicolumn{1}{|r|}{-0.4454} & \multicolumn{1}{|r|}{-0.008617} &
\multicolumn{1}{|r|}{-0.4540} & \multicolumn{1}{|r|}{3.54} &
\multicolumn{1}{|r|}{}\\\cline{1-7}%
0.2 & \multicolumn{1}{|r|}{-0.4442} & \multicolumn{1}{|r|}{5.62} &
\multicolumn{1}{|r|}{-0.4444} & \multicolumn{1}{|r|}{-0.01722} &
\multicolumn{1}{|r|}{-0.4616} & \multicolumn{1}{|r|}{1.93} &
\multicolumn{1}{|r|}{}\\\cline{1-7}%
0.3 & \multicolumn{1}{|r|}{-0.4369} & \multicolumn{1}{|r|}{7.17} &
\multicolumn{1}{|r|}{-0.4362} & \multicolumn{1}{|r|}{-0.02579} &
\multicolumn{1}{|r|}{-0.4620} & \multicolumn{1}{|r|}{1.85} &
\multicolumn{1}{|r|}{-0.4707}\\\cline{1-7}%
0.4 & \multicolumn{1}{|r|}{-0.4288} & \multicolumn{1}{|r|}{8.90} &
\multicolumn{1}{|r|}{-0.4278} & \multicolumn{1}{|r|}{-0.03433} &
\multicolumn{1}{|r|}{-0.4621} & \multicolumn{1}{|r|}{1.82} &
\multicolumn{1}{|r|}{}\\\cline{1-7}%
0.5 & \multicolumn{1}{|r|}{-0.4203} & \multicolumn{1}{|r|}{10.7} &
\multicolumn{1}{|r|}{-0.4202} & \multicolumn{1}{|r|}{-0.04280} &
\multicolumn{1}{|r|}{-0.4630} & \multicolumn{1}{|r|}{1.62} &
\multicolumn{1}{|r|}{}\\\hline
\end{tabular}
\end{center}
\label{tab_fourplates_sin}%
\end{table}

The relative errors versus the number of Gauss points and the WOS paths are
shown in Fig. \ref{fig_pathcomp}. The BIE-WOS result of $N_{g1}\times$
$N_{g1}=20\times20$, $N_{g2}\times$ $N_{g2}=10\times10$ , $N_{path}%
=2\times10^{3}$ and $a=0.5$ is taken as the reference solution. From Fig.
\ref{fig_pathcomp}, we can see that when the number of the Brownian paths
$N_{path}$ is larger than $10^{3}$ at each Gauss point, the BIE-WOS result
with $10\times10$ Gauss points will achieve an accuracy about 1\% in the
relative error.\qquad

\bigskip

\noindent$\bullet$ \textit{\bigskip Test 4: CPU time comparison with the
last-passage method}

\begin{figure}[ptb]
\begin{center}
\includegraphics[width=2.5in ]{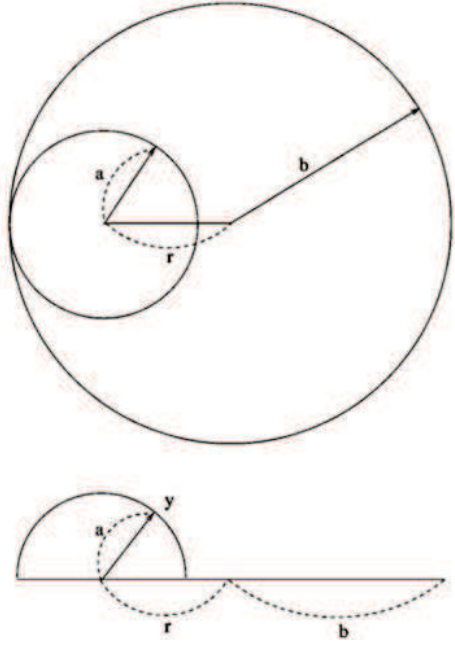}
\end{center}
\caption{Finding the charge distribution over a disk in 3-D}%
\label{fig_disk}%
\end{figure}

For both the last-passage and BIE-WOS methods, the CPU time is expected to be
linear in terms of the total number of random paths. We demonstrate this fact
with a case of a thin circular disk with radius $b$ in 3-D space
\cite{given02} as shown in Fig. \ref{fig_disk}. From \cite{given03}, the
analytical result of the charge density on the disk is:%

\begin{equation}
\sigma(\rho)=\frac{Q}{4\pi b\sqrt{b^{2}-\rho^{2}}},\ \ Q=8b. \label{bw95}%
\end{equation}

For a given relative error tolerance on the charge density at
$(-0.5,0,0)$, the CPU time comparison of both methods\ versus the number of
random paths are listed in Table \ref{tab_disk}. We take the radius $a=0.4$
for $S_{a},b=1$ for the radius of the thin disk, and the analytical charge
density is $\sigma(0.5)=0.735105$. From Table \ref{tab_disk}, we can see that
the CPU times are indeed in proportion to the total number of random paths for
both methods for a comparable accuracy. Though the integral $\Sigma
_{2}^{\prime}$ of the BIE-WOS method in this case is obviously zero, we still
evaluate it just as for a general variable potential and the CPU times of
$\Sigma_{2}^{\prime}$ is included in the CPU times of the BIE-WOS method in
Table \ref{tab_disk}. It is noted that the CPU times in computing the integral
$\Sigma_{2}^{\prime}$ for all cases are insignificant at about 0.012 second
for a 20$\times$20 Gauss quadrature.

\begin{table}[ptb]
\caption{The relative errors and the cpu times comparison according to the
number n of random paths}%
\label{tab_disk}
\begin{center}%
\begin{tabular}
[c]{|c|c|c|c|c|c|c|c|}\hline
\multicolumn{4}{|c}{Last-passage (LP)} & \multicolumn{4}{|c|}{BIE-WOS}\\\hline
$N_{\mathrm{path-LP}}$ & $\Sigma_{\text{LP}}$ & err$\%$ & cpu &
$N_{\mathrm{path-bie-wos}}$ & $\Sigma_{1}^{^{\prime}}+\Sigma_{2}^{^{\prime}}$
& err$\%$ & cpu\\
&  &  & time(s) &  &  &  & time(s)\\\hline
\multicolumn{1}{|r|}{$10^{4}$} & \multicolumn{1}{|r|}{0.69975} &
\multicolumn{1}{|r|}{-4.81} & \multicolumn{1}{|r|}{32} &
\multicolumn{1}{|r|}{$10^{2}\cdot100=10^{4}$} & \multicolumn{1}{|r|}{0.68888}
& \multicolumn{1}{|r|}{-6.29} & \multicolumn{1}{|r|}{30}\\\hline
\multicolumn{1}{|r|}{$10^{5}$} & \multicolumn{1}{|r|}{0.73253} &
\multicolumn{1}{|r|}{-0.35} & \multicolumn{1}{|r|}{331} &
\multicolumn{1}{|r|}{$10^{2}\cdot1000=10^{5}$} & \multicolumn{1}{|r|}{0.73960}
& \multicolumn{1}{|r|}{0.61} & \multicolumn{1}{|r|}{307}\\\hline
\multicolumn{1}{|r|}{$4\cdot10^{5}$} & \multicolumn{1}{|r|}{0.73743} &
\multicolumn{1}{|r|}{0.32} & \multicolumn{1}{|r|}{1325} &
\multicolumn{1}{|r|}{$20^{2}\cdot1000=4\cdot10^{5}$} &
\multicolumn{1}{|r|}{0.73441} & \multicolumn{1}{|r|}{-0.09} &
\multicolumn{1}{|r|}{1218}\\\hline
\end{tabular}
\end{center}
\end{table}

\subsection{Finding Neumann data over a patch of a curved boundary}

\begin{figure}[ptb]
\begin{center}
\includegraphics[width=2.5in ]{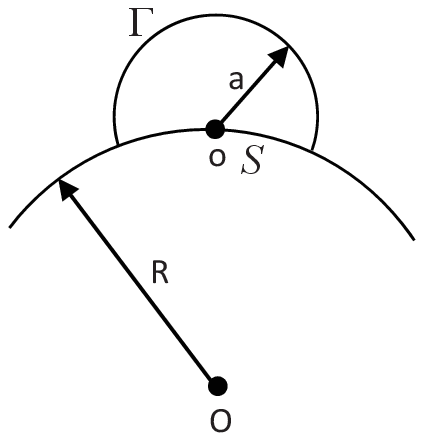}
\includegraphics[width=2.5in ]{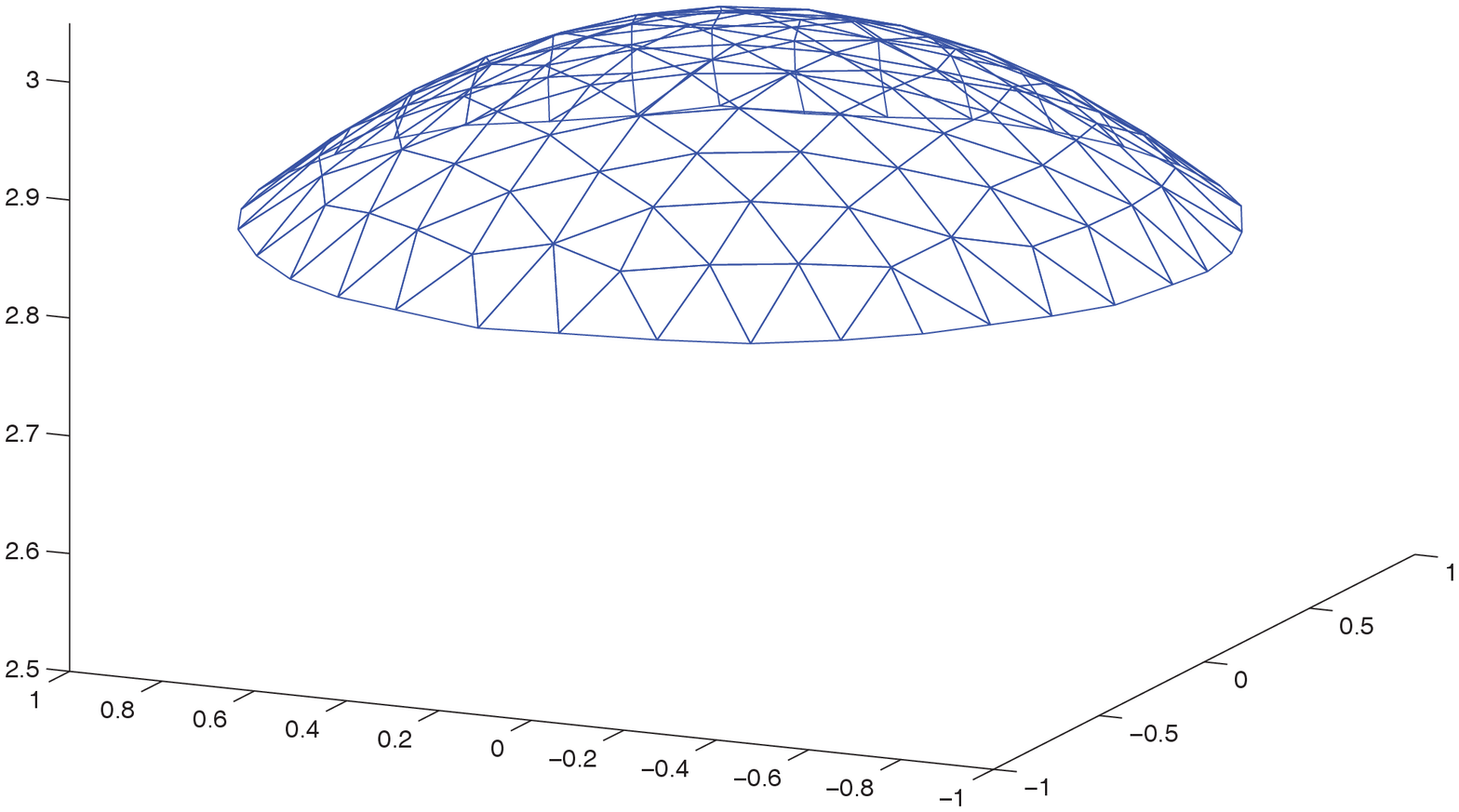}
\end{center}
\caption{Left: BIE-WOS setting for finding the Neumann data over a patch of a
big sphere; right: the mesh over the patch}%
\label{fig_bigsphere}%
\end{figure}

Next, to test the BIE-WOS method for a curved boundary, we compute the DtN
mapping on a big sphere as shown in Fig. \ref{fig_bigsphere} (left) with a radius $R=3$. A
point charge $q=1$ is located at the central point $O$ and the analytical
result for the potential is then known. To compute the Neumann data over a
local patch $S$ around the point $\mathbf{o}=(0,0,3)$ on the big spherical
surface, a small sphere with a radius $a=1$ is superimposed over the point
$\mathbf{o}$. The local patch $S$ is discretized with a triangular mesh as
shown in Fig. \ref{fig_bigsphere}.

The BIE equation of (\ref{bw49}) is solved by a collocation boundary element
method. When a collocation point is not inside an integration panel of the
BEM, a simple Gauss quadrature method is used. For collocation points inside
an integration panel, both weak and strong singularities will occur; however,
they can be regularized by a local polar transformation technique and a
$20\times20$ Gauss quadrature will then be used. For the integrals on $\Gamma
$, a $30\times30$ Gauss quadrature will be used. The potential $u(\mathbf{y})$
on $\Gamma$ is first computed, by the Feynman-Kac formula and the WOS\ method
with $10^{4}$ Brownian paths, on a regular grid, which is generated by evenly
discretizing the spherical surface along the polar and azimuthal angles. The
value $u(\mathbf{y})$ on $\Gamma$ but not on the grid points as required by
numerical quadratures will be interpolated using the values on the grid points.

\begin{figure}[ptb]
\begin{center}
\includegraphics[width=3.5in ]{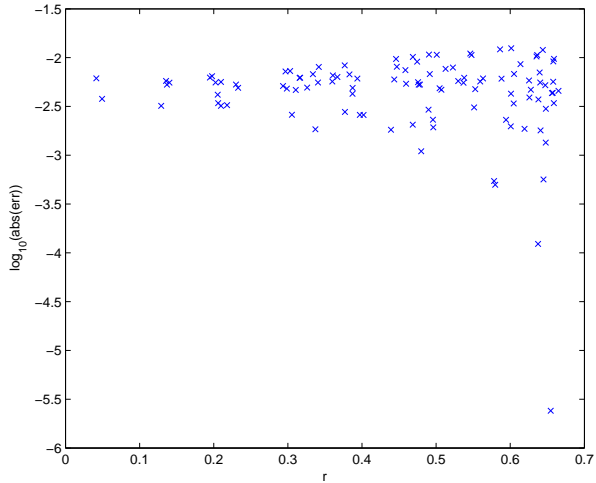}
\end{center}
\par
.\caption{Accuracy of the Neuamann data by the BIE-WOS method over the patch
$S,r<0.7a$, where $a$ is the radius of $\Gamma.$}%
\label{fig_accbigsphere}%
\end{figure}

The relative errors at the center of triangular panels on $S$ are shown in
Fig. \ref{fig_accbigsphere}, where the $x$-axis means the distance between the
center of the triangle to point $\mathbf{o}$. From Fig. \ref{fig_accbigsphere}%
, we can see that for the panels close to point $\mathbf{o}$, i.e. $r<0.7a$,
\ the maximal relative error is less than $1.25\%$, which will be accurate
enough for most engineering applications. It is noted that due to the sharp
corner edge singularity of the domain $\Omega_{S}$ where the hemisphere and
$\partial\Omega$ intersect, the piecewis constant collocation BEM will lose
some of its accuracy, which limits the region where acceptable accuracy of the
BEM solution can be used for the sought-after Neumann data. This well-known
problem in singular boundary elements usually is addressed with graded mesh
near the edge singularity \cite{chandler84}\cite{atkinson90}\cite{kress90} and
is still an active research topic in BEM methods \cite{bremer10}. A resolution
of this edge singularity can increase the region of useful BEM solution in the
BIE-WOS algorithm and can be incorporated into the algorithm. As discussed in
the last section, as the boundary $\partial\Omega$ will be covered with an
overlapping patch $S_{i}$, the loss of the accuracy of the BIE solution near
the edge of each patch will not hinder the use of the BIE-WOS method. However,
any improvement of the BEM near the edge will reduce the total number of
patches to cover the boundary, thus reducing the total cost.

\section{Conclusions and discussions}

\bigskip In this paper we have proposed a local BIE-WOS method which combines
a local deterministic singular BIE method and the Monte Carlo WOS algorithm to
find the Neumann data on general surfaces given Dirichlet data there. The
singular integral equation for the Neumann data at any single point or a local
patch on the boundary surface involves potential solution on a local
hemisphere, which can be readily obtained with the Feynman-Kac formula with
the help of the WOS sampling of the Brownian paths. Numerical results validate
the efficiency and accuracy of this method.

The local BIE-WOS method of finding the DtN or NtD mapping can give a parallel
algorithm for the solution of the Poisson equation with Dirichlet or Neumann
boundary conditions. Firstly, we partition the whole boundary $\partial\Omega$
into a union of overlapping patches $S_{i}$ namely,%

\[
\partial\Omega=\cup_{i}S_{i}.
\]
Then, the local BIE-WOS method can be used to find the DtN or NtD mapping over
each patch $S_{i}$ independently in parallel. In principle, the computation of
the BIE-WOS method over each patch can be done over one processor without need
for communications with others; thus a high parallel scalability can be
achieved. Secondly, the solution to the Poisson equation in the whole space
can be found with the integral representation of (\ref{bw03}) with the help of
one application of FMM \cite{Greengard}.

There are several important research issues to be addressed before the BIE-WOS
method can be used for large scale computation of Poisson or modified
Helmholtz equations. The first issue is (1) the NtD mapping problem, where the
Neumann data is given on the boundary and the Dirichlet data is required. In
this case, the Feynman-Kac formula derived in \cite{Hsu85} can be used, which
will involve reflecting Brownian paths \cite{Lions84} with respect to the
domain boundary. Efficient numerical implementation will have to be developed;
The second issue is the modified Helmholtz equation. Even though the
Feynman-Kac formula (\ref{bw11}) still applies, a survival factor will be
introduced as the WOS samples the Brownian paths and efficient ways to use the
Feynman-Kac formula will have to be addressed. The third issue is that since
the WOS scheme requires the computation of the distance between a Brownian
particle and the boundary of the solution domain, efficient algorithms will
have to be be studied for the overall speed of the BIE-WOS method.

The parallel algorithm based on the BIE-WOS method for solving Poisson or
modified Helmholtz equations will have the following important features:

\begin{itemize}
\item Non-iterative in construction and no need to solve any global linear system.

\item Stochastic in nature based on the fundamental link between the Brownian
motion and the solution of elliptic PDEs.

\item Massive parallelism suitable to large number of processors for large
scale computing due to the random walk and local integral equation components
of the algorithm.

\item No need for traditional finite element type surface or volume meshes.

\item Applicable to complex 3-D geometry with highly accurate treatment of
domain boundary geometries.
\end{itemize}

In comparison with traditional finite element and finite difference methods,
the BIE-WOS solver is only suitable for Poisson and modified Helmholtz
equations (due to the use of WOS-type sampling technique of the diffusion
paths) and its accuracy is limited to that of the Monte Carlo sampling
technique, while the traditional grid based methods can handle more general
PDEs with variable coefficients and achieve high accuracy. Nonetheless, as the
Poisson and modified Helmholtz equations form the bulk computation of
projection-type methods for incompressible flows and other important
scientific computing, the progress in scalability of parallel BIE-WOS
based-solvers will have large impact on the simulation capability of
incompressible flows and engineering applications.

\bigskip

\begin{center}
\textbf{Acknowledgement}
\end{center}

\medskip

The second author thanks the support of U.S. Army Research Office (grant
number W911NF-11-1-0364) and NSF (grant DMS-1005441). The first and third
authors are supported in part by NSFC research projects 60976034, 61076033,
61274032 and 61228401, the National Basic Research Program of China under the
grant 2011CB309701, the National major Science and Technology special project
2011ZX 01034-005-001-03 of China during the 12-th five-year plan period, the
Doctoral Program Foundation of the Ministry of Education of China
200802460068, the Program for Outstanding Academic Leader of Shanghai, and the
State Key Lab. of ASIC \& System Fudan Univ. research project 11MS013.

Finally, the authors acknowledge the constructive suggestions of the reviewers
on the presentation and organization of the materials, which have contributed
much to the improvement of the paper.

\bibliographystyle{siam}
\bibliography{bibfile}

\begin{thebibliography}{10}

\bibitem{atkinson90}
{\sc K.~Atkinson and I.~Graham}, {\em Iterative variants of the nystr\"{o}m
  method for second kind boundary integral operators}, SIAM J. Sci.Stat.
  Comput., 13 (1990), pp.~694--722.

\bibitem{brandt82}
{\sc A.~Brand}, {\em Guide to multigrid development}, Lecture Notes in
  Mathematics, 960 (1982), pp.~220--312.

\bibitem{brebbia78}
{\sc C.~A. Brebbia}, {\em The boundary element method in engineering}, Pentech
  Press. London, 1978.

\bibitem{bremer10}
{\sc J.~Bremer and V.~Rokhlin}, {\em Efficient discretization of laplace
  boundary integral equations on polygonal domain}, Journal of Computational
  Physics, 229 (2010), pp.~2507--2525.

\bibitem{cai13}
{\sc W.~Cai}, {\em Computational Methods for Electromagnetic Phenomena:
  Electrostatics in Solvation, Scattering, and Electron Transport}, Cambridge
  University Press, London, 2013.

\bibitem{chandler84}
{\sc GA~Chandler}, {\em Galerkin's method for boundary integral equations on
  polygonal domains}, The Journal of the Australian Mathematical Society.
  Series B. Applied Mathematics, 26 (1984), pp.~1--13.

\bibitem{zhu12}
{\sc G.~Chen, H.~Zhu, T.~Cui, Z.~Chen, X.~Zeng, and W.~Cai}, {\em Parafemcap: A
  parallel adaptive finite-element method for 3-d vlsi interconnect capacitance
  extraction}, Microwave Theory and Techniques, IEEE Transactions on, 60
  (2012), pp.~218--231.

\bibitem{chorin68}
{\sc A.J. Chorin}, {\em Numerical solution of the navier-stokes equations},
  Math. Comp, 22 (1968), pp.~745--762.

\bibitem{Chung95}
{\sc K.~L. Chung}, {\em Green, Brown and Probability}, world scientific, 1995.

\bibitem{Freidlin85}
{\sc M.I. Freidlin}, {\em Functional Integration and Partial Differential
  Equations.(AM-109)}, vol.~109, Princeton University Press, 1985.

\bibitem{friedman06}
{\sc A.~Friedman}, {\em Stochastic Differential Equations and Applications},
  Dover Publications, 2006.

\bibitem{nedelec78}
{\sc J.~Giroire and JC~Nedelec}, {\em Numerical solution of an exterior neumann
  problem using a double layer potential}, Math. Comp, 32 (1978), pp.~973--990.

\bibitem{given03}
{\sc J.A. Given and C.O. Hwang}, {\em Edge distribution method for solving
  elliptic boundary value problems with boundary singularities}, Physical
  Review E, 68 (2003), p.~046128.

\bibitem{given02}
{\sc J.A. Given, C.O. Hwang, and M.~Mascagni}, {\em First-and last-passage
  monte carlo algorithms for the charge density distribution on a conducting
  surface}, Physical Review E, 66 (2002), p.~056704.

\bibitem{Greengard}
{\sc L.~Greengard and V.~Rokhlin}, {\em A fast algorithm for particle
  simulations}, Journal of computational physics, 73 (1987), pp.~325--348.

\bibitem{Hsu85}
{\sc PEI HSU}, {\em Probabilistic approach to the neumann problem},
  Communications on pure and applied mathematics, 38 (1985), pp.~445--472.

\bibitem{hwang06}
{\sc C.O. Hwang and J.A. Given}, {\em Last-passage monte carlo algorithm for
  mutual capacitance}, Physical Review E, 74 (2006), p.~027701.

\bibitem{kress90}
{\sc R.~Kress}, {\em A nystr{\"o}m method for boundary integral equations in
  domains with corners}, Numerische Mathematik, 58 (1990), pp.~145--161.

\bibitem{coz98}
{\sc YL~Le~Coz, HJ~Greub, and RB~Iverson}, {\em Performance of random-walk
  capacitance extractors for ic interconnects: a numerical study}, Solid-State
  Electronics, 42 (1998), pp.~581--588.

\bibitem{coz92}
{\sc YL~Le~Coz and RB~Iverson}, {\em A stochastic algorithm for high speed
  capacitance extraction in integrated circuits}, Solid-State Electronics, 35
  (1992), pp.~1005--1012.

\bibitem{Lions84}
{\sc P.L. Lions and A.S. Sznitman}, {\em Stochastic differential equations with
  reflecting boundary conditions}, Communications on Pure and Applied
  Mathematics, 37 (1984), pp.~511--537.

\bibitem{Mascag03}
{\sc M.~Mascagni and C.O. Hwang}, {\em $\epsilon$-shell error analysis for
  ¡°walk on spheres¡± algorithms}, Mathematics and computers in simulation, 63
  (2003), pp.~93--104.

\bibitem{mascagni04}
{\sc M.~Mascagni and N.A. Simonov}, {\em The random walk on the boundary method
  for calculating capacitance}, Journal of Computational Physics, 195 (2004),
  pp.~465--473.

\bibitem{Muller56}
{\sc M.E. Muller}, {\em Some continuous monte carlo methods for the dirichlet
  problem}, The Annals of Mathematical Statistics, 27 (1956), pp.~569--589.

\bibitem{white91}
{\sc K.~Nabors and J.~White}, {\em Fastcap: A multipole accelerated 3-d
  capacitance extraction program}, Computer-Aided Design of Integrated Circuits
  and Systems, IEEE Transactions on, 10 (1991), pp.~1447--1459.

\bibitem{white94}
{\sc JR~Phillips and J.~White}, {\em A precorrected-fft method for capacitance
  extraction of complicated 3-d structures}, in Proceedings of the 1994
  IEEE/ACM international conference on Computer-aided design, IEEE Computer
  Society Press, 1994, pp.~268--271.

\bibitem{shebalfeld}
{\sc K.K. Sabelfeld and KK~Sabel}, {\em Monte Carlo methods in boundary value
  problems}, vol.~274, Springer-Verlag Berlin/Heidelberg/New York, 1991.

\bibitem{shebalfeld95}
{\sc KK~Sabelfeld, IA~Shalimova, and A.~Lavrentjeva}, {\em Random walk on
  spheres process for exterior dirichlet problem}, Monte Carlo Methods and
  Applications, 1 (1995), pp.~325--331.

\bibitem{shi02}
{\sc W.~Shi, J.~Liu, N.~Kakani, and T.~Yu}, {\em A fast hierarchical algorithm
  for three-dimensional capacitance extraction}, Computer-Aided Design of
  Integrated Circuits and Systems, IEEE Transactions on, 21 (2002),
  pp.~330--336.

\bibitem{temam84}
{\sc R.~Temam}, {\em Navier¨CStokes Equations: Theory and Numerical Analysis},
  North-Holland, Amsterdam, 1984.

\bibitem{widlund05}
{\sc A.~Toselli and O.~Widlund}, {\em Domain decomposition methods-algorithms
  and theory}, vol.~34, Springer, 2004.

\bibitem{shi05}
{\sc S.~Yan, V.~Sarin, and W.~Shi}, {\em Sparse transformations and
  preconditioners for 3-d capacitance extraction}, Computer-Aided Design of
  Integrated Circuits and Systems, IEEE Transactions on, 24 (2005),
  pp.~1420--1426.

\bibitem{yu04}
{\sc W.~Yu and Z.~Wang}, {\em Enhanced qmm-bem solver for three-dimensional
  multiple-dielectric capacitance extraction within the finite domain},
  Microwave Theory and Techniques, IEEE Transactions on, 52 (2004),
  pp.~560--566.

\end{thebibliography}

\end{document}